\documentclass{elsarticle}

\usepackage{graphicx}
\graphicspath{{figures/}}
\usepackage[fleqn]{amsmath} % fleqn eqs a la  izquierda
\usepackage{amsfonts}
\usepackage{parskip}
\usepackage{color}
\usepackage{dsfont}
\usepackage{siunitx} % para manejar consistentemente notacion de numeros y unidades:
\usepackage{booktabs} % para tablas mas bonitas
\usepackage{placeins} % para \FloatBarrier
\usepackage{subcaption} % para subfiguras y subcaptions
\usepackage[nonumberlist, nopostdot, nogroupskip, acronym]{glossaries}  % para administrar simbolos y acronimos
        \setglossarystyle{super}  % estilo de glosario (2 columnas justificado izquierda)
        \newglossary*{symbol}{Symbols}  % nuevo glosario para symbolos
\newacronym{pod}{POD}{proper orthogonal decomposition}
\newacronym{pgd}{PGD}{proper generalized decomposition}
\newacronym{fe}{FE}{finite element}
\newacronym{fe2}{FE$^2$}{finite element square}
\newacronym{svd}{SVD}{singular value decomposition}
\newacronym{hprfe2}{HPR-FE$^2$}{high-performance reduced finite element square}
\newacronym{rve}{RVE}{representative volume element}
\newacronym{rom}{ROM}{reduced order model}
\newacronym{roec}{REOC}{reduced energy-based optimal cubature}
\newacronym{hf}{HF}{high-fidelity finite element}
\newacronym{hffe2}{HF-FE$^2$}{high-fidelity finite element square}
\newacronym{ibvp}{IBVP}{initial boundary value problem}
\newacronym{bvp}{BVP}{boundary value problem}
        \makenoidxglossaries

\usepackage[normalem]{ulem}

\begin{document}
\begin{frontmatter}
\title{HIGH PERFORMANCE REDUCTION TECHNIQUE FOR MULTISCALE FINITE ELEMENT MODELING (HPR-FE$^{2}$): TOWARDS INDUSTRIAL MULTISCALE FE SOFTWARE}
%
% Group authors per affiliation:
\author[a]{Marcelo Raschi}
\author[a,c]{Oriol Lloberas-Valls}
\author[b,c]{Alfredo Huespe}
\author[a,c]{Javier Oliver}

\address[a]{CIMNE -- Centre Internacional de Metodes Numerics en Enginyeria, Campus Nord UPC, Mòdul C-1, c/ Jordi Girona 1-3, 08034, Barcelona, Spain}
\address[b]{CIMEC-UNL-CONICET, Predio Conicet, Ruta Nac. 168 s/n - Paraje El Pozo, 3000, Santa Fe, Argentina}
\address[c]{E.T.S d'Enginyers de Camins, Canals i Ports, Technical University of Catalonia (BarcelonaTech), Campus Nord UPC, M\`odul C-1, c/ Jordi Girona 1-3, 08034, Barcelona, Spain}


\begin{abstract}
The authors have shown in previous contributions that reduced order modeling with optimal cubature applied to \gls{fe2} techniques results in a reliable and affordable multiscale approach, the \gls{hprfe2} technique.
Such technique is assessed here for an industrial case study of a generic 3D reinforced composite whose microstructure is represented by two general microcells accounting for different deformation mechanisms, microstrucural phases and geometry arrangement.
Specifically, in this approach the microstrain modes used for building the \gls{rom} are obtained through standard \gls{pod} techniques applied over snapshots of a representative sampling strain space.
Additionally, a reduced number of integration points is obtained by exactly integrating the main free energy modes resulting from the sampling energy snapshots.
The outcome consists of a number of dominant strain modes integrated over a remarkably reduced number of integration points which provide the support to evaluate the constitutive behavior of the microstructural phases.
It is emphasized that stresses are computed according to the selected constitutive law at the reduced integration points and, therefore, the strategy inherits advantageous properties such as model completeness and customization of material properties.
Overall results are discussed in terms of the consistency of the multiscale analysis, customization of the microscopic material parameters and speedup ratios compared to \gls{hf} simulations.
\end{abstract}

\end{frontmatter}

%%%%%%%%%%%%%%%%%%%
% Section symbols and acronyms
%
\glsresetall
 % reset the use status of all acronyms after Abstract
\clearpage
%\printnoidxglossary[type=symbol]
\printnoidxglossary[type=acronym]
%
%%%%%%%%%%%%%%%%%%%

%%%%%%%%%%%%%%%%%%%
% Section
\section{Introduction}
Modern material design industry demand high-performing simulation tools capable of capturing the non-linear effects of materials observed at the structural scale, although governed by the physics occurring at small length scales.
During the last years, a considerable effort has been put into the development and improvement of multiscale computational homogenization techniques for analyzing heterogeneous composite materials~\cite{matouvs2017review}.

Typically, for material modeling involving two well-separated length scales, e.g., structural and mesoscopic scales, special attention has been paid to hierarchical techniques that show a high capability to capture the complex mechanical interaction effects between scales (\cite{michel1999effective}, \cite{Miehe_et_al_1999}, \cite{terada2001class}, \cite{blanco2016variational}).
A pioneering paper of this kind of techniques is the \gls{fe2} approach by Feyel et al.~\cite{feyel_fe2_2000}, which can tackle the non-linear behavior of complex microstructures by solving surrogate problems on a \gls{rve} of the target material~\cite{Gitman_et_al_2007,Nguyen_et_al_2010}.

However, \gls{fe2} approaches are still computationally unaffordable for their use in software for practical applications.
This is the main reason why they have not yet been widely transferred to the industrial sector.
As an example, an analysis involving a standard \gls{fe2} multiscale simulations may typically take from months to years  to complete in typical computing clusters.

A general concept to defeat the barrier imposed by the overwhelming cost of \gls{fe2} techniques, without giving up the displayed accuracy to capture the interaction effects between scales, consists of separating the full computational burden of a multiscale material simulation in two stages.
First, an offline training stage is carried out with the \gls{hf} model of a \gls{rve}, where a data set of \gls{hf} solutions is sampled and a surrogate or reduced model of the homogenized material response is built.
Finally, this reduced model, which demands a low computational cost, is used in the online structural virtual testing.
Several methodologies that follow this concept are very briefly described:

\begin{itemize}
\item [{\it i)}] Artificial neural networks can be used to construct a cheap parameterized surrogate model.
Results reported in~\cite{Rocha_2020} clearly show that the computational cost of these strategies during the online stage can be lower than the one demanded by alternative reduced models.
However, surrogate models derived in this way ignore the physical basis of the problem, as the model parameters are non-physical coefficients.
Additionally, these strategies are unable to simulate material problems governed by input parameters unforeseen in the sampled load trajectories;
e.g., for reproducing non-monotonous loading conditions that were not encountered in the sampling.
See~\cite{ghavamian2019accelerating}, where a technique to overcome this drawback in history-dependent materials is proposed.

\item [{\it ii)}] Employment of wavelet functions (\cite{Deslauriers_1989}, \cite{van_tuijl_2019}) to define different resolutions for the interpolation of basis functions.
This strategy can be seen as an on-the-fly optimal integration technique for reduced models, at the expense of a low speedup---only one order of magnitude---with respect to a full \gls{rom}.

\item [{\it iii)}] The \gls{pgd} method~\cite{Chinesta_2013} is a powerful reduction technique that has been used in a number of engineering applications.
Within the \gls{pgd} approach, a multidimensional solution is computed offline, where model parameters and boundary conditions are considered extra coordinates of the problem.
The multidimensional solution is approximated by the sum of modes, composed by the product of the dimensional coordinates (including space-time, material properties and boundary conditions).
This product consists of a lightweight offline computational catalog to be employed for rapid simulations of industrial interest.
Typical procedures adopting this approach have been addressed by Ladeveze and coworkers~\cite{ladeveze2010latin, chinesta2011short}, Doblare and coworkers~\cite{Doblare_2016}, among others.
\item [{\it iv)}] Hyper-reduced finite element methods, resulting from a combination of a \gls{rom} technique (which in turn is based on \gls{pod}) with a procedure that reduces the computational complexity related to the non-linear terms of the variational formulation.
The \gls{pod} technique provides an empirical low-dimensional basis of spatial interpolation functions that replaces the original \gls{fe} interpolation functions of the \gls{hf} model.
Furthermore, several techniques to decrease computational complexity have emerged in last years.
The objective is to select only a few Gauss integration points---the \emph{cubature} points---where the non-linear terms are computed.
The parameters of these cubature points are associated with a physical term of the problem.
Some examples of hyper-reduction techniques are the discrete empirical interpolation~\cite{chaturantabut2010nonlinear}, the energy conserving sampling and weighting method~\cite{zahr2017multilevel}, the use of reduced integration domain~\cite{ryckelynck2009hyper}, and the use of \gls{pod} modes to approximate the integrand of the non-affine terms, plus the selection of a reduced number of sampling points via gappy data reconstruction~\cite{hernandez2014high}.
A related technique to the present approach is the empirical cubature method reported in~\cite{J.A.Hernandez2016}, and improved in \cite{hernandez2020multiscale}.
\end{itemize}
Following the approach explained in point \emph{iv}, the \gls{hprfe2} technique---described by the authors in previous contributions \cite{OliverROM2017},~\cite{caicedo2019high}---is evaluated in this paper.
The principal characteristic of this model lies in its potential use as a software engine for the material design industry to overcome the \emph{tyranny of scales}.
It has been reported in~\cite{Lloberas_Complas_2019} that this \gls{hprfe2} technique is able to provide a speedup of up to four orders of magnitude for microstructural discretizations of several millions of integration points, assuming relative errors of around \SI{1}{\percent} with respect to \gls{hf} solutions.
These speedups imply that a multiscale simulations that originally could take up to several years, can now be performed in a few hours.
This clearly indicates that \gls{hprfe2} technology fulfills the requirements of the material design industry for an affordable multiscale analysis tool.

It is emphasized that the \gls{hprfe2} formulation proposed in this paper differs from the Empirical Cubature Method proposed in \cite{J.A.Hernandez2016} in the sense that, in order to construct a reduced integration scheme, what is approximated is the fundamental (variational) principle, i.e., the free energy minimization, and not its derivatives (stresses and internal forces) as it will be described in the formulation section.
The minimum of the variational principle (energy minimization) is, in any case, well approximated and since is a scalar field, the reduction procedure yields a lighter reduced basis.
Consequently, the performance of the approach is considerably improved with respect to the approach from~\cite{J.A.Hernandez2016}, which is based on integration of vector entities.
This is, among others, one of the specific features of this method.

The objective of this work is to further assess the potentiality of this technique in terms of a trade-off between computational speedup, consistency, accuracy, and preservation of the physical basis of the model.
The main novelty of the current contribution consists in the application of a HPR-FE$^2$ methodology to an industrial simulation case (fibre-reinforced composite laminates).
This entails a number of upgrades to the strategy with respect to former contributions: 3D analysis, tailor-made constitutive modeling, full reconstruction of stress and internal variables fields, assessment of the real speedups and computation times for realistic simulation scenarios.

In the following Sections \ref{Formulacion} and \ref{RFE2}, a brief summary of the \gls{hprfe2} technique is presented (a detailed description can be found in \cite{OliverROM2017}, \cite{caicedo2019high}, \cite{Lloberas_Complas_2019}).
The accuracy and consistency of the methodology is assessed in Section \ref{sect:results}, by comparing the solutions of two composites virtually tested with the \gls{hprfe2} and the \gls{hf} techniques.
%%%%%%%%%%%%%%%%%%%

%%%%%%%%%%%%%%%%%%%
% Section
%%%%%%%%%%%%%%%%%%%%%%%%%%%%%%
\section{Multiscale problem formulation}
\label{Formulacion}
\begin{figure}[htbp]
\centering
\includegraphics[width=0.8\linewidth]{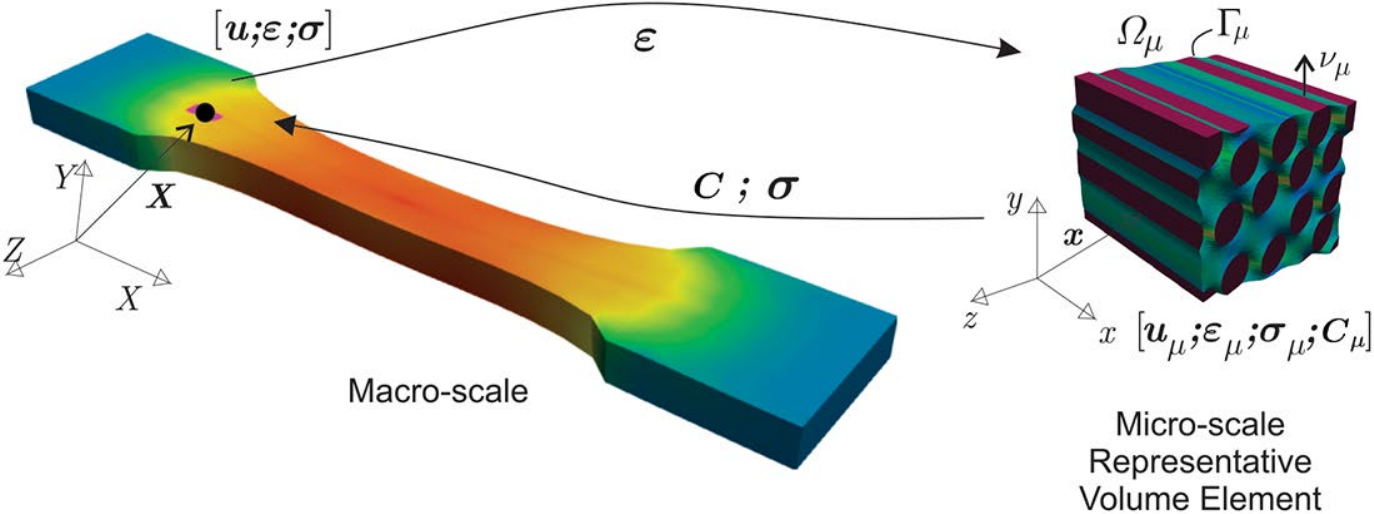}
\caption{\gls{fe2} multiscale material model.}
\label{fig_multiscale}
\end{figure}

%%%%%%%%%%%%%%%%%%%%%%%%%%%%%%
%
Let a hierarchical multiscale material model such as the one schematized in Figure~\ref{fig_multiscale} be assumed.
The coupon represents the macroscale structure, where effective mechanical variables are considered, i.e., displacement $\boldsymbol{u}$, strain $\boldsymbol{\varepsilon}$, and stress $\boldsymbol{\sigma}$.
A \gls{rve} of the microstructure is used to compute the homogenized constitutive response of the material at each point $\boldsymbol{X}$ of the structure.
Computational homogenization consists in transferring $\boldsymbol{\varepsilon}$ onto the \gls{rve}, solving a microscale \gls{bvp}, and upscale the effective stress $\boldsymbol{\sigma}$ evaluated as the volumetric average of the microstress $\boldsymbol{\sigma}_{\mu}$.
The same operation also computes the corresponding effective constitutive tangent tensor $\boldsymbol{C}$.

The first point to be defined in the formulation is the admissible kinematics at the microscale.
The displacement field at the microscale $\boldsymbol{u}_{\mu}$, results in the addition of three terms:
\begin{equation}
\boldsymbol{u}_{\mu}(\boldsymbol{x})
=
\boldsymbol{u} + \boldsymbol{\varepsilon} \cdot \boldsymbol{x} + \tilde{\boldsymbol{u}}_{\mu}(\boldsymbol{x}),
\label{displ_micro_0}
\end{equation}
where the last term is the displacement fluctuation field $\tilde{\boldsymbol{u}}_{\mu}(\boldsymbol{x})$.
Accordingly, the microstrain $\boldsymbol{\varepsilon}_\mu$ in the \gls{rve} is given by the addition of two terms:
\begin{equation}
    \boldsymbol{\varepsilon}_\mu (\boldsymbol{x})
     =
    \boldsymbol{\varepsilon} + \tilde{\boldsymbol{\varepsilon}}_{\mu}(\boldsymbol{x})
    =
    \boldsymbol{\varepsilon} + \nabla^{s}_{\boldsymbol{x}}\tilde{\boldsymbol{u}}_{\mu}(\boldsymbol{x})
    \quad
    \forall \, \boldsymbol{x} \in \Omega_\mu,
\label{epsilon_mu_0}
\end{equation}
where $\boldsymbol{\varepsilon}$ is assumed uniform in the microscale volume $\Omega_\mu$, while the microstrain fluctuation $\tilde{\boldsymbol{\varepsilon}}_{\mu}$ is the symmetric gradient of the displacement fluctuation field\footnote{In the remaining of this paper, the infinitesimal strain theory will be considered.}, satisfying:
\begin{equation}
    \int_{\Omega_\mu}\tilde{\boldsymbol{\varepsilon}}_{\mu} \; d\Omega =
    \int_{\Omega_\mu}\nabla^{s}_{\boldsymbol{x}}\tilde{\boldsymbol{u}}_{\mu} \; d \Omega =
    \int_{\Gamma_{\mu}} \tilde{\boldsymbol{u}}_{\mu} \otimes^{s} \boldsymbol{\nu}_{\mu} \; d \Gamma
    = \boldsymbol{0}.
\label{strain_fluc_zero_0}
\end{equation}

In accordance with this kinematics, a \gls{bvp} at the \gls{rve} is next formulated in terms of microstrain fluctuations~\cite{OliverROM2017}.

MICROSCALE PROBLEM FORMULATION:
Given the macroscale strains $\boldsymbol{\varepsilon}$, and the spaces of kinematically compatible strain fluctuations $\mathcal{U}_{\mu}^{\tilde{\varepsilon}}$ and admissible strain fluctuations $\mathcal{V}_{\mu}^{\tilde{\varepsilon}}:$
\begin{equation}
\mathcal{U}_{\mu}^{\tilde{\varepsilon}}
=
\mathcal{V}_{\mu}^{\tilde{\varepsilon}}
:=
\left \{
\tilde{\boldsymbol{\varepsilon}}_{\mu}
\quad|\quad
\int_{\Omega_\mu}\tilde{\boldsymbol{\varepsilon}}_{\mu}
\; d \Omega = \boldsymbol{0}
\quad \text{and} \quad
\tilde{\boldsymbol{\varepsilon}}_{\mu} \in \mathcal{E}_{\mu}
\right \} ;
\label{equilibrium_PI_1-1}
\end{equation}
find
$\tilde{\boldsymbol{\varepsilon}}_{\mu} \in \mathcal{U}_{\mu}^{\tilde{\varepsilon}}$
such that
\begin{equation}
\int_{\Omega_\mu}
\boldsymbol{\sigma}_{\mu} (\boldsymbol{\varepsilon}_\mu, d_{\mu}) :
\mathbf{\mathbf{\widehat{\boldsymbol{\varepsilon}}_{\mu}}}
\; d \Omega = 0  ; \qquad
\forall \; \mathbf{\widehat{\boldsymbol{\varepsilon}}_{\mu}} \in \mathcal{V}_{\mu}^{\tilde{\varepsilon}} ;
\label{equilibrium_PI_2-1}
\end{equation}
\begin{equation}
\dot{d}_{\mu}(\boldsymbol{x}, \boldsymbol{\varepsilon}_\mu)
=
g(\boldsymbol{\varepsilon}_\mu, d_{\mu}),
\label{damage evolution_PI_3-1}
\end{equation}
where the space $\mathcal{E}_{\mu}$ of microstrain tensorial functions is defined such that its elements fulfil the infinitesimal strain compatibility conditions given by
\begin{equation}
    \mathcal{E}_{\mu}:=
    \left\{
     \boldsymbol{\zeta} \in \mathbb{S}^{n_\text{dim}\times n_\text{dim}}
    \quad \left| \right. \quad
    e_{mjq} e_{nir} \boldsymbol{\zeta}_{ij, qr} =
    0
    \right \},
\label{compat_eq}
\end{equation}
with $n_\text{dim}=3$ being the dimension of the Euclidean space.

Equation~\eqref{equilibrium_PI_2-1} is the variational form of the equilibrium equation in terms of the microstress $\boldsymbol{\sigma}_{\mu}$ which, in turn, implicitly considers the constitutive relation as a function of the internal variable\footnote{In this work, a single scalar internal variable \emph{damage} is considered.} microscopic damage ${d}_{\mu}$ and microstrains.
Conventional damage constitutive models are assumed at the microscale.
Therefore, the damage rate equation is explicitly defined by \eqref{damage evolution_PI_3-1}.

The homogenized stress $\boldsymbol{\sigma}$ is the volumetric average of the microstresses which are the solution of the \gls{bvp} \eqref{equilibrium_PI_1-1}-\eqref{damage evolution_PI_3-1}:
\begin{equation}
\boldsymbol{\sigma}
=
\frac{1}{|\Omega_\mu|}\int_{\Omega_\mu} \boldsymbol{\sigma}_{\mu}(\boldsymbol{\varepsilon}_\mu, d_{\mu} )
\; d\Omega,
\end{equation}
and the homogenized constitutive tensor is:
\begin{equation}
\boldsymbol{C}(\boldsymbol{X}, t)
=
\frac{\partial\boldsymbol{\sigma}}{\partial\boldsymbol{\varepsilon}}
=
\frac{1}{|\Omega_\mu|}\int_{\Omega_\mu}\boldsymbol{C}_{\mu}(\boldsymbol{x})
\Big(\mathbb{I}+\boldsymbol{A}_{\mu}(\boldsymbol{X}, \boldsymbol{x}, t)\Big)
\; d{\Omega},
\label{C_effect}
\end{equation}
where $\boldsymbol{C}_{\mu}$ are the constitutive tensor of the composite phases, $\mathbb{I}$ is the fourth order identity tensor and $\boldsymbol{A}_{\mu}$ is the localization tensor ($\tilde{\boldsymbol{\varepsilon}}_{\mu} = \boldsymbol{A} \boldsymbol{\varepsilon}$).

Since the microscale \gls{bvp} is formulated in terms of the strain field, the displacement $\boldsymbol{u}$ does not play any role in the \gls{rve} problem, and the history of the strain components $\boldsymbol{\varepsilon}$ are the only parameters determining the homogenized response\footnote{This is also a feature of the proposed formulation, that allows optimal specific treatment of zones of the \gls{rve} with different constitutive equations.}.

By adopting \eqref{epsilon_mu_0} and a kinematically admissible microdisplacement fluctuation space satisfying the last identity in \eqref{strain_fluc_zero_0}, the micromechanical \gls{bvp} can be re-written in a more conventional displacement formulation.
As will be shown in the next section, both formulations, in displacements and strains, are used to develop the \gls{hprfe2} technique.
The displacement-based formulation is used in the \gls{hf} \gls{fe} technique for the model sampling stage, and the strain-based formulation is employed within the reduced and hyper-reduced models.
The advantage of taking this approach is commented in the following section.
%
%%%%%%%%%%%%%%%%%%%%%%%
\section{\Acrfull{hprfe2}}
\label{RFE2}
%%%%%%%%%%%%%%%%%%%%%%%
%
The \gls{hprfe2} model uses a low-dimensional space of functions for approximating the microfluctuation strain field $\tilde{\boldsymbol{\varepsilon}}_{\mu}$, as described in Section~\ref{sect:reduction}.
An identical approach is taken for the microstrain variations, i.e., a Galerkin formulation.
Additionally, a reduced numerical cubature rule, the \acrfull{roec}, is introduced in Section~\ref{Sec_ROEC} to compute the integral balance in~\eqref{equilibrium_PI_2-1}.
Section~\ref{sect:recovery} describes the procedure for recovering  the strains, displacements and damage fields of the \gls{rve}, i.e., projecting fields from the cubature points onto the original \gls{fe} mesh.
\subsection{Dimensional reduction of the microstrain fluctuation field}
\label{sect:reduction}
The low-dimensional space representing $\tilde{\boldsymbol{\varepsilon}}_{\mu}$  is built as the span of an orthogonal basis of $n_{\varepsilon}$ spatial functions with global support:
 $
 \{\boldsymbol{\Psi} (\boldsymbol{x})\}
 =
 \{ \Psi_1(\boldsymbol{x}), \dotsc, \Psi_{n_{\varepsilon}}(\boldsymbol{x})\}
 $
(where $\Psi_{i}(\boldsymbol{x}) \in \mathbb{R}^{n_{\sigma}}$ and $n_{\sigma}=6$ is the dimension of the strain and stress tensors in Voigt notation for 3D problems), as follows:
\begin{equation}
\tilde{\boldsymbol{\varepsilon}}_\mu(\boldsymbol{x},t)
=
\sum_{i=1}^{n_\varepsilon}\Psi_i(\boldsymbol{x})c_i(t)
=
\boldsymbol{\Psi}(\boldsymbol{x})\boldsymbol{c}(t),
\label{epsilon_reducido}
\end{equation}
where each element $\Psi_i$ of the basis $\left \{\boldsymbol{\Psi}\right \}$ is a microstrain fluctuation mode, and the vector of pseudo-time dependent coefficients $\boldsymbol{c}(t) = [c_1, \dotsc, c_{n_\varepsilon}]$
($\boldsymbol{c} \in \mathbb{R}^{n_\varepsilon}$) represents the amplitude of these modes.
In the last identity of \eqref{epsilon_reducido}, the matrix
$
\boldsymbol{\Psi}(\boldsymbol{x})
=
[\Psi_1, \dotsc, \Psi_{n_\varepsilon}]
$,
with $\boldsymbol{\Psi}(\boldsymbol{x})\in\mathbb{R}^{n_{\sigma}\times n_{\varepsilon}}$ collects, in columns, the $n_{\varepsilon}$ microstrain modes of the basis $\left\{\boldsymbol{\Psi}\right\}$.
Notice that, in order to preserve a simplified notation, an identical symbol $\tilde{(\boldsymbol{\cdot})}$ identifies the microstrain fluctuation field in both,  \gls{hf} and low-dimensional approaches.

The basis $\left\{\boldsymbol{\Psi}\right\}$, which is composed of $n_\varepsilon$ basis vectors (or modes), is computed with a \gls{pod} technique applied to a set of microstrain fluctuation fields, in turn obtained as solutions of the microcell problem (\gls{hf} model formulated in displacements) during an offline sampling process.

Finally, the variational formulation \eqref{equilibrium_PI_2-1}-\eqref{damage evolution_PI_3-1} is projected onto the space spanned by the basis $\left\{ \boldsymbol{\Psi}\right\}$.
Therefore, the number of equations of the corresponding discrete non-linear systems reduces to $n_\varepsilon$ equations, whose solution determines the vector $\boldsymbol{c}(t)$.

Remarkably, by construction, each element $\Psi_i$ belongs to the vectorial space  \eqref{compat_eq} inheriting the boundary conditions imposed to the \gls{hf} model in the sampling stage.
Therefore, any function spanned by the basis $\{\boldsymbol{\Psi}\}$ is an admissible microfluctuation strain.
%
%%%%%%%%%%%%%%%%%%%%%%%%%
\subsection{ \Acrfull{roec}}
\label{Sec_ROEC}
%%%%%%%%%%%%%%%%%%%%%%%%%%%%%%%%%%%%%%%%%%%%
%
The integral term in the global balance equation \eqref{equilibrium_PI_2-1} is computed with the \gls{roec}  rule.
This rule is derived as follows (find additional details in \cite{OliverROM2017}): the dimensionality of the  free energy function space $\phi_\mu$ computed in the sampling stage is reduced using a similar approach to that adopted for the microstrain fluctuation field \eqref{epsilon_reducido}.
Thus, assuming that the low-dimensional free energy is spanned by a basis of $N_\phi$ elements,
\begin{equation}
\phi_\mu(\boldsymbol{x}, t)
=
\sum_{i=1}^{N_\phi}\Phi_i (\boldsymbol{x}) c^\phi_i(t)
=
\boldsymbol{\Phi}(\boldsymbol{x})\boldsymbol{c}^\phi(t).
\label{phi_reducido}
\end{equation}
The number of cubature points of the \gls{roec} rule $N_r$, their spatial position $\boldsymbol{z}$ and the corresponding weights $\omega$ are selected with a similar criterion to that reported in \cite{J.A.Hernandez2016}.
This criterion is based on adopting an exact integration\footnote{Exact integration understood as the selection of the optimal set (among the Gauss integration rule adopted in the \gls{hf} model) and its corresponding weights, in order to exactly integrate the selected $N_r = N_\phi + 1$ energy modes $\boldsymbol{\Phi}$.}
of  the free energy basis $\Phi_{i}$, plus the condition that the reduced integration of the unit-function in $\Omega_\mu$ gives the volume of the microcell $\sum_{j=1}^{N_r} \omega_j = | \Omega_\mu |$.
These conditions are expressed as follows:
\begin{eqnarray}
\int_{\Omega_\mu}^\text{Gauss}\Phi_i \; d\Omega
& \simeq &
\int_{\Omega_\mu}^\text{\gls{roec}} \Phi_i \; d\Omega
=
\sum_{j=1}^{N_r}\Phi_i(z_j) \; \omega_j
\qquad \forall \; i = 1, \dotsc, N_\phi,
\label{eq:reduced_quadrature} \\
\int_{\Omega_\mu}^\text{Gauss} \; d\Omega
& \simeq&
  \int_{\Omega_\mu}^\text{\gls{roec}} \; d\Omega
  =
 \sum_{j=1}^{N_r} \; \omega_j =| \Omega_\mu |,
\end{eqnarray}
where the symbol $\int_{\Omega_\mu}^\text{Gauss} (\cdot) \; d\Omega $ refers to the Gauss quadrature rule in the \gls{hf} model.
This criterion provides $N_\phi + 1$ conditions.
Therefore, the number of cubature points which can be obtained is $Nr = N_\phi + 1$.

\noindent {\bf Remark:} A remarkable feature of the \gls{roec} technique hinges on the determination of the reduced integration rule in terms of the free energy field (fundamental variational principle) and not on its derivatives (stresses and internal forces), as proposed in many other formulations.
As a result, the problem goal, i.e., find the minimum of the variational principle (free energy), is the actual key of the solution of the problem, unlike in alternative methods based on finding null values of the functional derivatives (internal forces)\footnote{I.e., the Euler-Lagrange equations of the variational principle.}.
Consequently, the formulation results simpler and very accurate, leading to the high performance reduction technique \gls{hprfe2}.
%
%%%%%%%%%%%%%%%%%%%%%%%
\subsection{Reconstruction of fields in the \gls{rve}}
\label{sect:recovery}
%%%%%%%%%%%%%%%%%%%%%%%
%
After solving the \gls{hprfe2} problem, given the approximate microstrain field  $\tilde{\boldsymbol{\varepsilon}}_{\mu}$ and the internal variables  of the damage model ($r^\text{\gls{roec}}$ described in \ref{sect:appendix}) in the cubature points, the microfields
$(\tilde{\boldsymbol{u}}_{\mu}, \; \boldsymbol{\sigma}_{\mu}, \; d_{\mu})$
can be recovered and projected onto the original \gls{hf} mesh as described below.
\subsubsection{Reconstruction of the displacement fluctuation field }
\label{sect:recovery_displacement}
The displacement fluctuation fields $\tilde{\boldsymbol{u}}_{\mu}$ can be recovered at the original nodes of the \gls{hf} mesh by using the following procedure.
Let the conventional space of the \gls{fe} interpolation functions for displacement $\mathcal{V}^\text{HF}$ associated to the \gls{hf} mesh be defined, but excluding the displacements of rigid body modes.
Then, the displacement reconstruction problem is stated as follows:

FIND: $\tilde{\boldsymbol{u}}_{\mu} \in \mathcal{V}^\text{HF}$ such that:
\begin{equation}
\int_{\Omega_\mu}^\text{Gauss}
(\underbrace{\tilde{\boldsymbol{\varepsilon}}_{\mu}}_{\boldsymbol{\Psi}\boldsymbol{c}}
-
\nabla^{s}\tilde{\boldsymbol{u}}_{\mu})
:
\nabla^{s}\hat{\boldsymbol{u}}_{\mu} \; d\Omega = 0;  \quad
\forall \hat{\boldsymbol{u}}_{\mu} \in \mathcal{V}^\text{HF}.
\label{ppp}
\end{equation}
Introducing the \gls{fe} approach of $\tilde{\boldsymbol{u}}_{\mu}$ and the  variations of microdisplacement fluctuations $\hat{\boldsymbol{u}}_{\mu}$,
and considering that $\tilde{\boldsymbol{q}}_{\mu}$ is the global vector of nodal parameters, interpolating $\tilde{\boldsymbol{u}}_{\mu}$ in the original \gls{fe} mesh, then \eqref{ppp} can be re-written as follows:
\begin{equation}
\underbrace{
    \left(
        \int_{\Omega_\mu}^\text{Gauss} \boldsymbol{B}_{u}^{T}\boldsymbol{\Psi} \, d \Omega
    \right)
}_{\mathbb{F}_u}\boldsymbol{c}
-
\underbrace{
    \left(
        \int_{\Omega_\mu}^\text{Gauss}\boldsymbol{B}_{u}^{T}\boldsymbol{B}_{u}d\Omega
    \right)
}_{\mathbb{K}_u}
\tilde{\boldsymbol{q}}_{\mu}
 = 0,
\end{equation}
where $\boldsymbol{B}_{u}$ is the strain-displacement matrix of the conventional \gls{hf} \gls{fe} method ($\nabla^{s}\tilde{\boldsymbol{u}}_{\mu}=  \boldsymbol{B}_{u}\tilde{\boldsymbol{q}}_{\mu}$).
Finally:
\begin{equation}
\tilde{\boldsymbol{q}}_{\mu}
=
\underbrace{{\mathbb{K}_u}^{-1}\mathbb{F}_u}_{\mathbb{D}_{\boldsymbol{u}}}\boldsymbol{c}
=
\mathbb{D}_{\boldsymbol{u}}\boldsymbol{c}.
\label{last_eq}
\end{equation}

Notice that matrix $\mathbb{D}_{\boldsymbol{u}}$ from \eqref{last_eq} has to be computed only once for the reconstruction process.
\subsubsection{Reconstruction of the internal variable field of the damage model}
Similarly, the internal variable $r$ of the damage model (summarized in \ref{sect:appendix}) can also be recovered.
Firstly, a low-dimension spatial approach for $r$ is defined.
In the offline sampling of the original microcell, snapshots of this field are gathered and a reduced basis $\boldsymbol{\Psi}^r$ of $n_r$ elements is determined by means of a \gls{pod} technique, similarly to the procedure defined in \eqref{epsilon_reducido}.
Then, the reduced field $r$ can be expressed as:
\begin{equation}
r(\boldsymbol{x} ,t)
=
\sum_{i=1}^{n_r} \Psi^r_{i}(\boldsymbol{x}) c^r_{i}(t)
=
\boldsymbol{\Psi}^r(\boldsymbol{x}) \boldsymbol{c}^r(t).
\label{r_reducido}
\end{equation}

The vector $\boldsymbol{r}^\text{\gls{roec}} \in \mathbb{R}^{N_\text{r}}$ (whose components are the internal variable values at the $N_\text{r}$ reduced cubature points and at a given pseudo-time $t$) is obtained from the \gls{hprfe2} solution.
Secondly, $\boldsymbol{r}^\text{\gls{roec}}$ is projected onto the vector $\boldsymbol{r}^\text{\gls{hf}} \in \mathbb{R}^{N_\text{gp}}$ (being $N\text{gp}$ the total number of Gauss integration points), whose components are the internal variables at the Gauss quadrature points of the \gls{hf} mesh.
The following problem is defined:

%\vskip 0.3truecm
FIND:
$\boldsymbol{c}^r\in \mathbb{R}^{n_r}$ such that:
\begin{equation}
\int_{\Omega_\mu}^\text{Gauss}
(\boldsymbol{\Psi}^r \, \boldsymbol{c}^r
-
\boldsymbol{r}^\text{\gls{roec}}):\boldsymbol{\Psi}^r\, \delta \boldsymbol{c}^r \; d \Omega = 0
\quad
\forall \delta \boldsymbol{c}^r \in \mathbb{R}^{n_r}.
\label{int_var}
\end{equation}

Applying a similar procedure to that used for the microdisplacement fluctuation recovery in Section~\ref{sect:recovery_displacement}, coefficients $\boldsymbol{c}^r$ can be computed as follows:
\begin{equation}
\boldsymbol{c}^r = \underbrace{{\mathbb{K}_r}^{-1}\mathbb{F}_r}_{\mathbb{D}_{r}} \boldsymbol{r}^\text{\gls{roec}}
 =
 \mathbb{D}_{r} \boldsymbol{r}^\text{\gls{roec}}; \quad
\mathbb{F}_r= \int_{\Omega_\mu}^\text{Gauss}(\boldsymbol{\Psi}^r)^{T}\, d\Omega; \quad
\mathbb{K}_r= \int_{\Omega_\mu}^\text{Gauss}(\boldsymbol{\Psi}^r)^{T}\boldsymbol{\Psi}^r\, d\Omega.
\label{last_eqrvi}
\end{equation}

Finally, using equation \eqref{r_reducido}, we recover the vector $\boldsymbol{r}$, whose $i$-th component $r_i$ (at the spatial position of the $i$-th Gauss point $\boldsymbol{x}_i$) satisfies
\begin{equation}
r_i(\tau_i) =
\begin{cases}
0                             &  \quad \text{if} \quad \tau_i \leq 0;  \\
\tau_i                    & \quad \text{if} \quad  0 < \tau_i < 1;  \\
1                             &  \quad\text{if} \quad   1 \leq \tau_i ,   \\
\end{cases}
\label{r_reducido_2}
\end{equation}
where $\tau_i = \boldsymbol{\Psi}^r(\boldsymbol{x}_i)\boldsymbol{c}^r$.
Since the damage variable is recovered by projecting a vector from a lower dimensional approximation space onto the higher dimensional space of finite elements, in this operation there may arise values lesser than 0 or greater than 1.
Thus, it is necessary to disregard those values assigning to the recovered damage variable the value 0 or 1, respectively, to each case.
\subsubsection{Reconstruction of the microstress tensor field}
Given the internal variable vector $\boldsymbol{r}$  and the microstrain
$\boldsymbol{\varepsilon}_\mu = \boldsymbol{\varepsilon} + \boldsymbol{\Psi}\boldsymbol{c} $
at each Gauss quadrature point of the \gls{hf} \gls{fe} mesh, we compute the microstress tensor $\boldsymbol{\sigma}_{\mu}$ at the same points by appealing to the damage constitutive equation.

\noindent {\bf Remark:} Most of the involved computations in \eqref{last_eq} and \eqref{last_eqrvi} can be performed during the offline model sampling stage, typically, the evaluation of matrices $\mathbb{D}_{\boldsymbol{u}}$ and $\mathbb{D}_{r}$.
The matrix vector products
$\mathbb{D}_{\boldsymbol{u}} \boldsymbol{c}$ and $\mathbb{D}_{r} \boldsymbol{r}^\text{\gls{roec}}$,
as well as  the evaluation of the microstress at the Gauss integration points of the \gls{hf} \gls{fe} mesh can be computed during a post-processing stage  and after solving the \gls{hprfe2} problem.
These additional evaluations represent a low computational cost.
%%%%%%%%%%%%%%%%%%%

%%%%%%%%%%%%%%%%%%%
% Section
\section{\Acrfull{hprfe2} model assessment}
\label{sect:results}
The \gls{hprfe2} technique is assessed by means of the multiscale simulation using two microstructures of engineering interest, representing glass fiber-reinforced composites with epoxy matrix.
Particular attention is paid to the trade-off between attained accuracy, computational speedup, and the capacity of the \gls{hprfe2} model for simulating non-sampled trajectories.

The offline strategy of the sampling stage of both models is explained in Section \ref{sect:training}, which details the number and directions of the loading scenarios needed to attain the reduced basis employed in different approximations of the \gls{hf} solution.
The error of the reduced order approximation with respect to the \gls{hf} solution is calculated as explained in Section~\ref{sect:errors}.
\emph{Model completeness}, i.e., the capacity for modeling material responses which have not been specifically considered during the sampling stage, is detailed in Section~\ref{sect:completeness}.
\emph{Computational performance} for microcells of increasing size and complexity is shown in Section~\ref{sect:speedup}.
\emph{Material customization}, i.e., the capacity of the model to reproduce the response of materials whose properties are different to those used during the offline model sampling, is evaluated in Section~\ref{sect:customization}.

Examples of multiscale virtual tests are shown in Section~\ref{sect:dogbone}, where a test coupon is simulated under varying configurations of the microstructure (obtained by rotations of the same microcell and customization of the material parameters), and a number of combinations of representation modes and integration points.
Recovered damage and stress fields of the microcell in sampled elements of the coupon is also shown.

%
%%%%%%%%%%%%%%%%%%%%%%%%%%%%%%%%%%%%%%%%%%%%
\subsection{Microscopic models}
\label{sect:models}
%%%%%%%%%%%%%%%%%%%%%%%%%%%%%%%%%%%%%%%%%%%%
%
%
\begin{itemize}

\item [{\it i)}] {\bf {Model \emph{A}}} represents a periodic composite constituted by a generic cross-ply laminate containing several plies of aligned longitudinal fibers (\SI{40}{\micro\meter} diameter, volume fraction of \SI{35}{\percent}) at \ang{0} and \ang{90} angles, embedded in an epoxy matrix.
It represents a material that can manifest the modeled failure mechanisms (fiber-matrix pull-out, fiber-matrix decohesion, and interply delamination), rather than a real-world material example.
The \gls{rve} of this composite, as well as the material properties of its constituent phases, are defined in Table~\ref{tab:material_params_ref}.
Fibers are modeled elastic, while matrix, interply, and fiber-matrix cohesive interphase are assumed inelastic, modeled by a damage constitutive model with bilinear hardening, as described in \ref{sect:appendix}.

\item [{\it ii)}] {\bf Model \emph{B}} represents one ply of an industrial multilayer composite.
The microstructure consists of a random distribution of unidirectional fibers (\SI{7}{\micro\meter} diameter, volume fraction of \SI{60}{\percent}) in an epoxy matrix.
Fibers are assumed elastic in the strain range considered.
Matrix is modeled by the same damage model referred previously.
Contrarily to model \emph{A}, the fiber-matrix interaction effects ignore the fiber pull-out and decohesion mechanisms.
The material parameters of the constituent phases of this composite are reported in Table~\ref{tab:material_params_porto}.
\end{itemize}

The objective pursued with the analysis of both models is not necessarily to perform their precise assessment and validation against experimental results.
Rather, it is the demonstration of the capabilities of the \gls{hprfe2} technique for capturing the composite deformation modes involving discontinuous strain fields, particularly those arising in model \emph{A}.
Accurately capturing these deformation modes results in a challenging target for the present approach.
To assess the \gls{hprfe2} capability, the solutions obtained with this methodology are compared to those obtained with the \gls{hf} model, which are considered as the reference solutions in all tested cases.

\begin{table}[hbt]
    \centering
    \begin{tabular}{lcccc}
        \toprule
        \centering{
            \includegraphics[height=2.5cm]{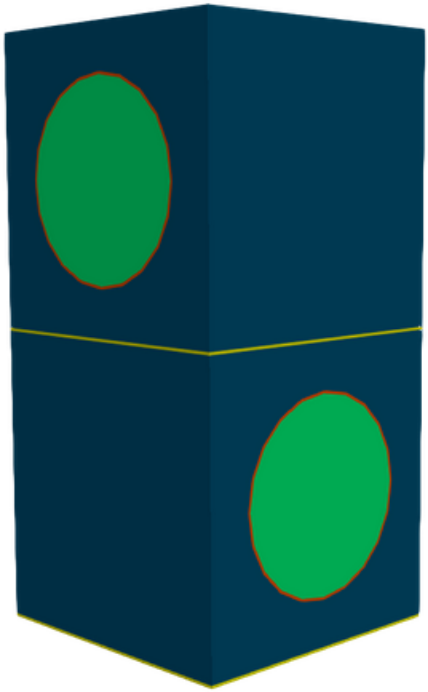}
            \includegraphics[width=.07\textwidth]{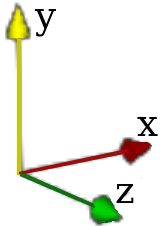}
        }    &
        \includegraphics[height=2.5cm]{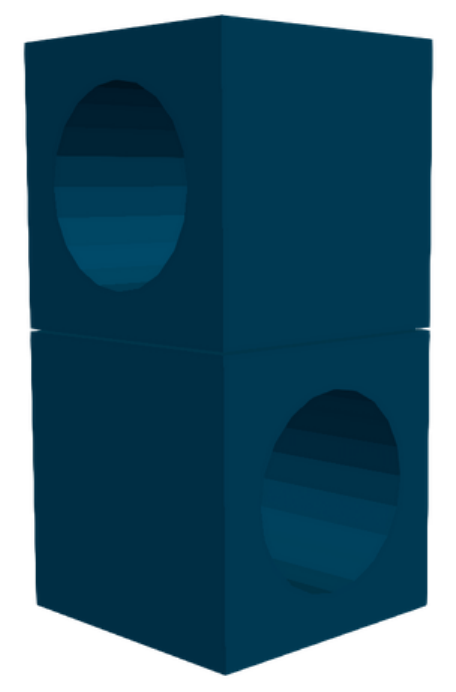}        &
        \includegraphics[height=2.5cm]{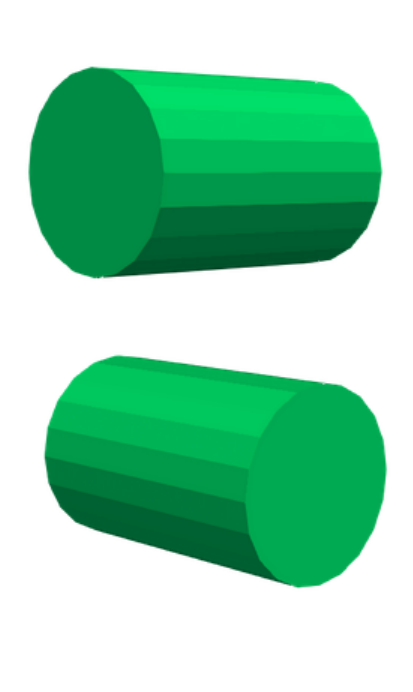}        &
        \includegraphics[height=2.5cm]{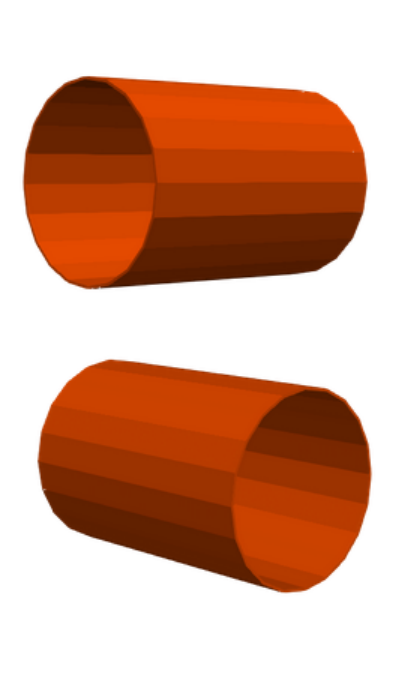}        &
        \includegraphics[height=2.5cm]{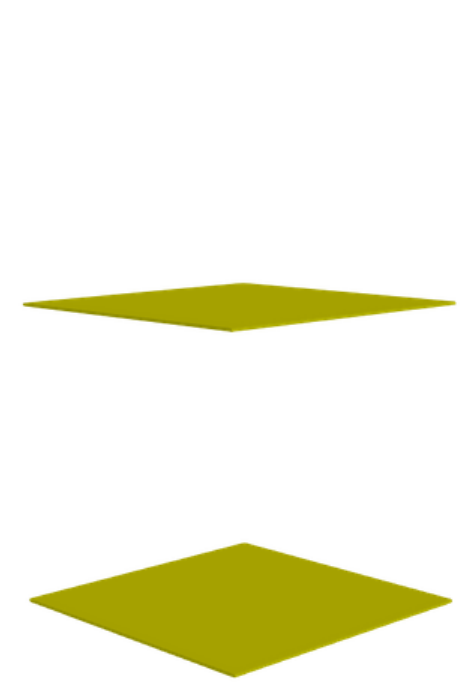}        \\
        & Matrix     & Fiber & Cohesive    & Interply    \\
        \midrule
        Young modulus $E$ [\si{N/mm^2}] & \num{3790} & \num{72000} & \num{3790}  & \num{3790}  \\
        Poisson ratio $\nu$& \num{0.37} & \num{0.25}  & \num{0.37}  & \num{0.37}  \\
        Elastic threshold $\sigma_{\text{e}}$ [\si{N/mm^2}] & \num{200}  & -- & \num{100}   & \num{100}   \\
        Infinity elastic threshold $\sigma_\infty$ [\si{N/mm^2}] & \num{220}  & -- & \num{100.1} & \num{100.1} \\
        Hardening parameters $H_1$, $H_2$& \num{0.10}, \num{0.01} & --  & \num{0.01}, \num{0.01}  & \num{0.01}, \num{0.01}  \\
        \bottomrule
    \end{tabular}
    \caption{{\bf Model \emph{A}.}
\gls{rve}.
Composite phases: matrix, fiber (volume fraction \SI{35}{\percent}), fiber-matrix cohesive (\SI{0.8}{\percent}), interply (\SI{0.5}{\percent}).
Material parameters of the constituent phases according to damage model described in \ref{sect:appendix}.
Failure mechanisms modeled: damage in matrix regions, delamination between plies, and fiber-matrix pull-out and debonding.
Fibers are assumed elastic.
}
    \label{tab:material_params_ref}
\end{table}

\begin{table}
    \centering
    \begin{tabular}{lcc}
        \toprule
        \centering{
            \includegraphics[height=2.5cm]{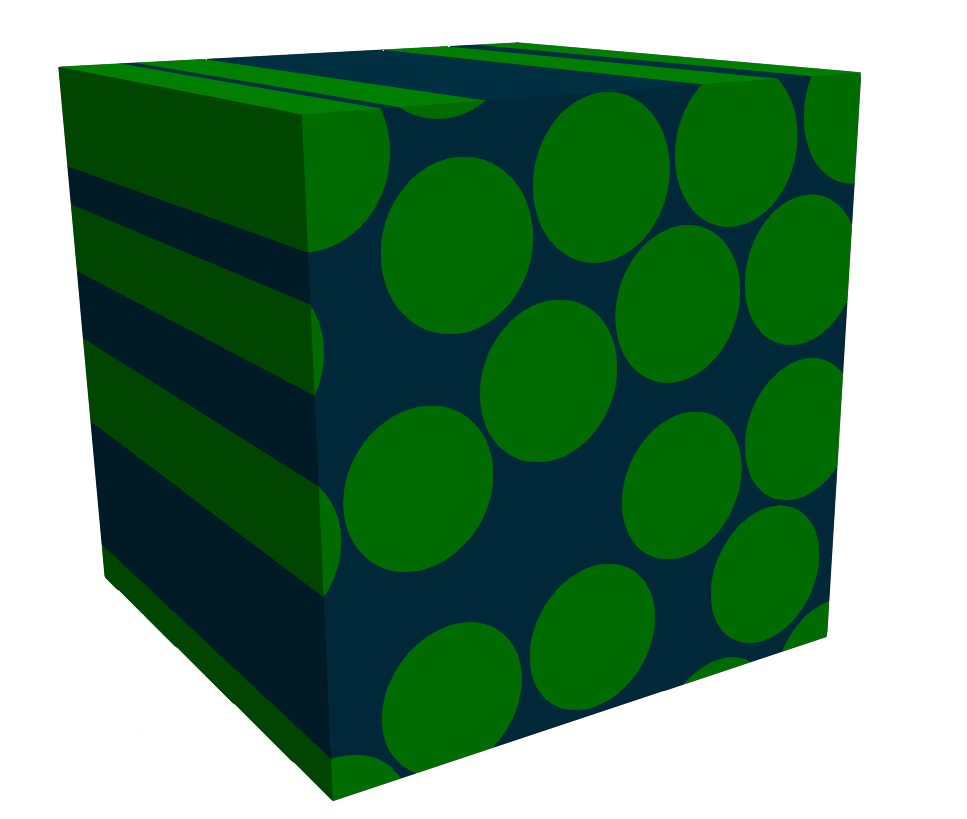}
            \includegraphics[width=.07\textwidth]{ref_local_4}
        }    &
        \includegraphics[height=2.5cm]{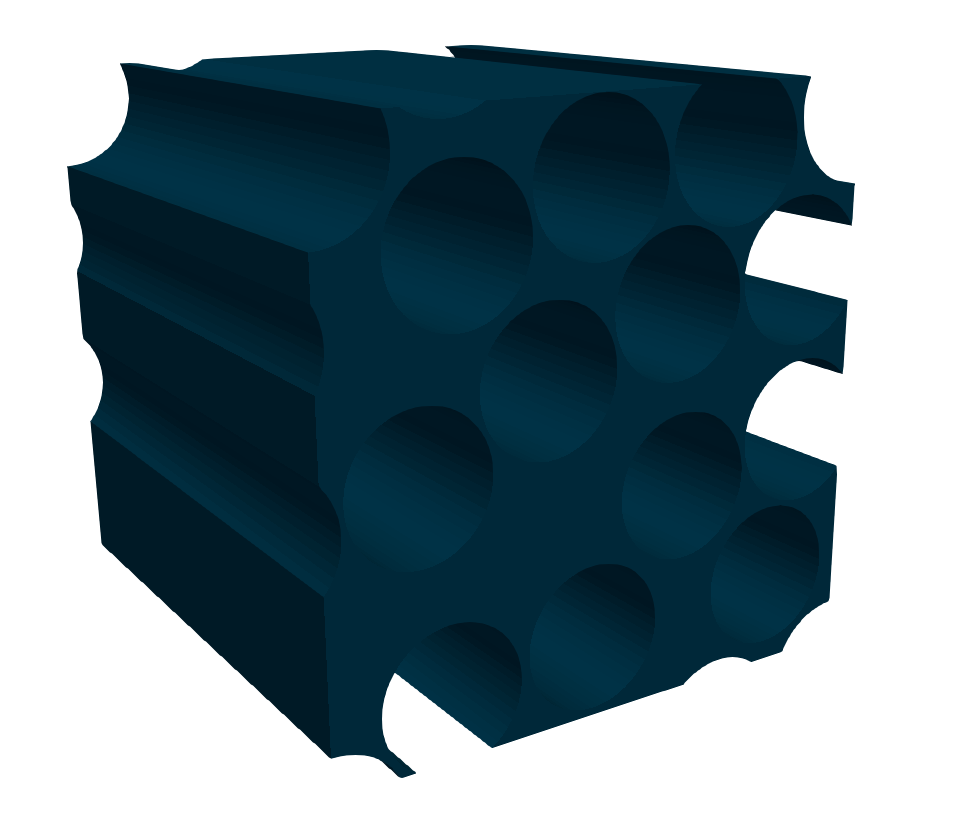}        &
        \includegraphics[height=2.5cm]{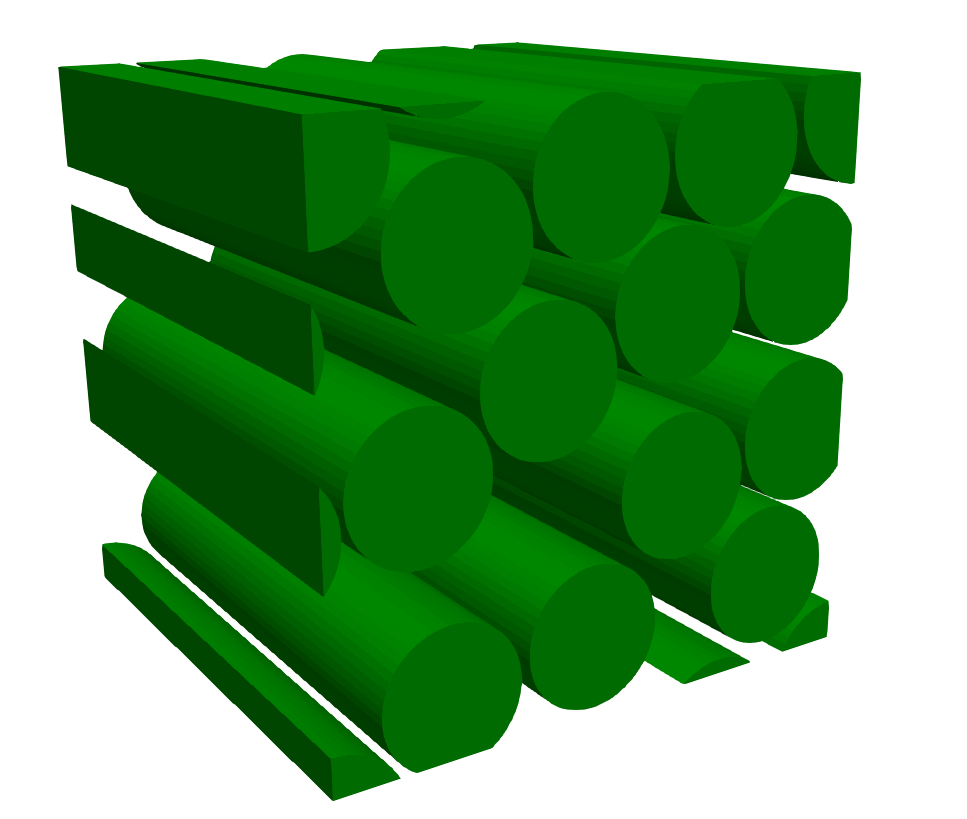}        \\
          & Matrix & Fiber \\
        \midrule
        Young modulus $E$ [\si{N/mm^2}]                                        & \num{4000} & \num{231000} \\
        Poisson ratio $\nu$                                                                      & \num{0.38}  & \num{0.2}  \\
        Elastic threshold $\sigma_{\text{e}}$ [\si{N/mm^2}]       & \num{60}  & -- \\
        Infinity elastic threshold $\sigma_\infty$ [\si{N/mm^2}] & \num{70}  & -- \\
        Hardening parameter $H_1$, $H_2$                                                    & \num{0.335}, \num{0.05}  & -- \\
        \bottomrule
    \end{tabular}
    \caption{{\bf Model \emph{B}.} Microcell. Composite phases: matrix, fiber (\SI{60}{\percent}). Material parameters of the constituent phases according to damage model described in \ref{sect:appendix}. Fibers are assumed elastic. In all the microcells of these models, the boundary material topology is enforced to be periodic.}
    \label{tab:material_params_porto}
\end{table}
\FloatBarrier
%%%%%%%%%%%%%%%%%%%%%%%%%%%%%%%%%%%%%%%%%%%%
\subsection{Sampling }
\label{sect:training}
%%%%%%%%%%%%%%%%%%%%%%%%%%%%%%%%%%%%%%%%%%%%
%
The sampling stage is performed by solving a set of \emph{trajectories}, i.e., problems in the \gls{rve} with the \gls{hf} technique, subjected to a macroscopic strain as the evolving action.

By considering that the six components of a symmetric macrostrain tensor can be related to a point lying on the $\mathbb{R}^6$ space, we choose \num{100} points uniformly distributed in the unitary hypersphere in $\mathbb{R}^6$ \cite{fritzen2019}, to define the same number of corresponding macrostrains.
Macrostrains are imposed to the \gls{rve} multiplied by a time factor $\chi$, monotonically increasing from 0 to $\chi_\text{end} = 0.1$ for model \emph{A}, and $\chi_\text{end} = 0.02$ for model \emph{B}.
The coefficient $\chi_\text{end}$ is chosen such that most of the trajectories reach the inelastic regime, which assures a sufficient number of inelastic snapshots during the sampling stage.
The solutions of the full loading trajectory set constitutes the \gls{rve} sampling procedure.
Although a damage model is employed to capture the material non-linearities, only a moderate strain localization is expected since softening is not included in the constitutive response.
 Hence, the adopted sampling range suffices to correctly reproduce the non-linear regime and it is not expected that severe localization takes place way beyond the one captured during the sampling stage.
  An extension of the method to account for softening would require a regularization strategy such as the one analyzed in\cite{OliverROM2017}.

The sampling data set is composed of \num{40} snapshots (one for each timestep) per each trajectory, from where the \gls{pod} bases are obtained.
There are no special requirements for the selection of the snapshots, provided that the selected sampling trajectory represents the full parametric strain space.
The approach taken here consists in gathering as many inelastic snapshots as possible, in a large trajectory diversity, and to leave the selection of optimal snapshots to the SVD stage.
A snapshot consists of strain, energy and damage variables in every Gauss integration point (ip) of the \gls{hf} mesh, at a given time step.
These snapshots are required for building the reduced set of integration (\emph{cubature}) points and the \gls{pod} strain modes for the \gls{hprfe2} model, and the reconstruction of the displacement and damage fields of the \gls{rve}.

Snapshots in the elastic regime are processed separately from those obtained during the inelastic regime.
%
%\rem{Details of this procedure can be seen in \cite{OliverROM2017} and \cite{caicedo2019high}}\nota{(la referencia a la teoria ya se hace al final de la intro, pag 4)}.
%
This strategy guarantees an exact response of the \gls{hprfe2} model in the elastic regime.

The sampling was performed using KratosMultiphysics simulation software \cite{kratos}.
\FloatBarrier
%%%%%%%%%%%%%%%%%%%%%%%%%%%%%%%%%%%%%%%%%%%%
\subsection{Accuracy: assessing non-sampled trajectories}
\label{sect:errors}
%%%%%%%%%%%%%%%%%%%%%%%%%%%%%%%%%%%%%%%%%%%%
%
The accuracy of the model is primarily assessed by comparing, for a particular non-sampled validation trajectory, the homogenized stresses obtained with the reduced model and the \gls{hf} simulation.

The \gls{bvp} of the microcells of models \emph{A} and \emph{B} are solved with the \gls{hprfe2} technique using a combination of modes in the range of \SIrange{10}{70}{}, and cubature points in the range of \SIrange{100}{2600}{}.
The obtained results correspond to an imposed monotonic macrostrain trajectory, defined by
\begin{equation}
\boldsymbol{\varepsilon}(\chi) = [
\num{-0.076}, \num{0.748}, \num{0.188}, \num{0.539}, \num{0.006}, \num{-0.329}
]\, \chi,
\label{eq:validation_trajectory}
\end{equation}
until reaching the time $\chi=\chi_\text{end}$.
This trajectory is not included in the sampling process.

\begin{figure}[hbt]
    \centering
    \includegraphics[width=\textwidth]{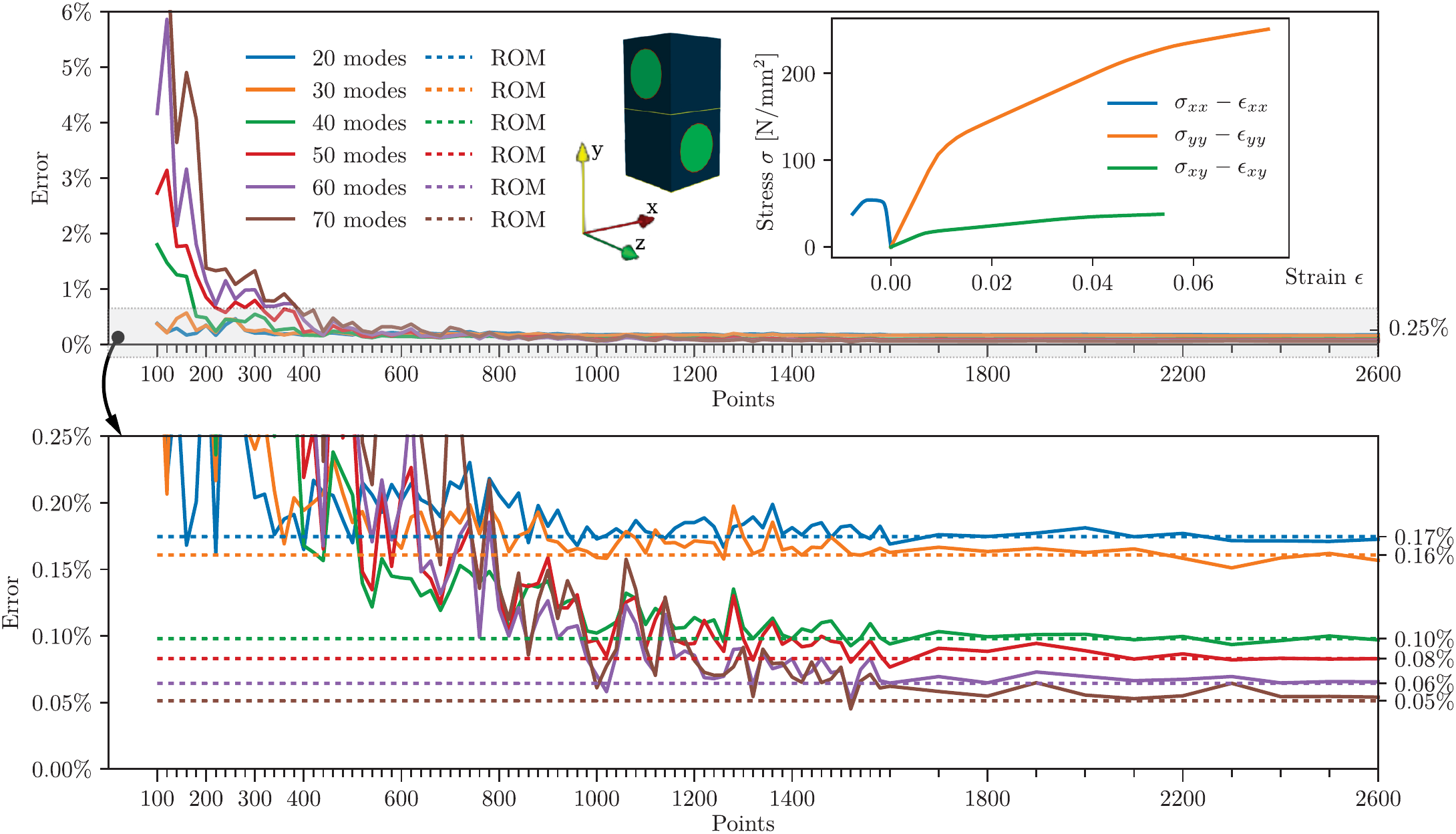}
    \caption{{\bf Model \emph{A}.} Error versus number of cubature points of the \gls{roec} scheme, for different numbers of strain modes. The dotted lines correspond to \gls{rom} errors (with exact \gls{hf} Gauss integration rule). The bottom figure is a detailed view of the top figure.}
    \label{fig_errors_ref}
\end{figure}

\begin{figure}[hbt]
    \centering
    \includegraphics[width=\textwidth]{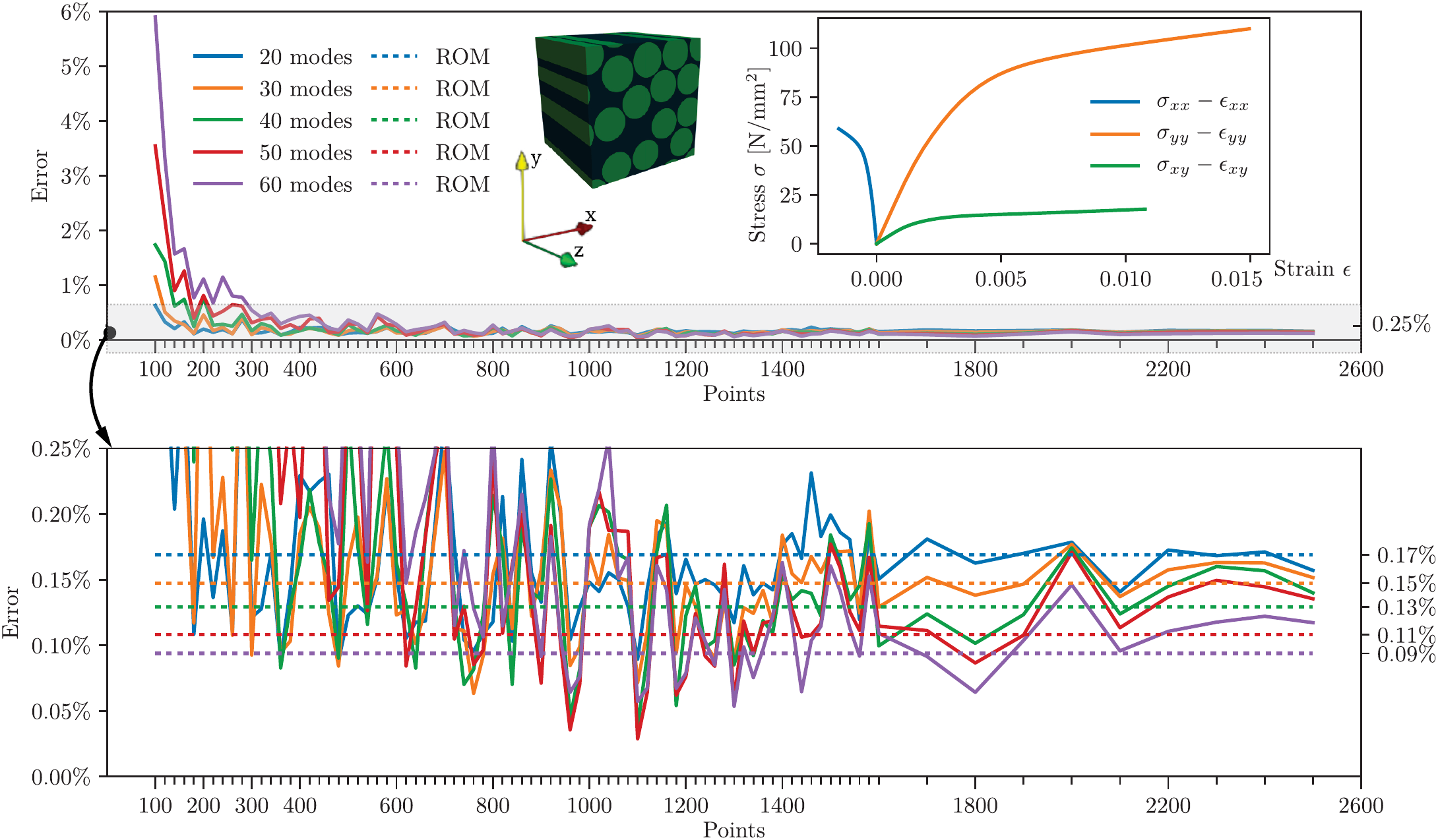}
    \caption{{\bf Model \emph{B}.} Error versus number of cubature points of the \gls{roec} schemes, for different numbers of strain modes. The dotted lines correspond to ROM errors (with exact \gls{hf} Gauss integration rule).  The bottom figure is a detailed view of the top figure.}
    \label{fig_errors_porto}
\end{figure}

The error for each of these combinations is defined as
\begin{equation}\label{eq_error}
\text{Error}
=
\max_{\chi} \frac
{ \| \sigma_\text{HF}(\chi) - \sigma_\text{R}(\chi) \|_{\infty} }
{ \|\sigma_\text{HF}(\chi)\|_{\infty}}
 \qquad \text{for} \; \chi\in [0,\chi_\text{end}],
\end{equation}
i.e., the $L^{\infty}$ norm of the maximum relative difference between the homogenized stress tensor given by the reduced model $\sigma_\text{R}$, and the \gls{hf} model stress $\sigma_\text{HF}$, computed along the trajectory.
Errors, given in terms of the number of cubature points and strain modes, are plotted in Figures~\ref{fig_errors_ref} and~\ref{fig_errors_porto} for the  microcells of models \emph{A} and \emph{B}.
In general, and as expected, we can observe that a higher number of strain modes reduces the approximation error provided that a high enough number of cubature points is used.
Also, as expected, the curves tend to an asymptotic value coinciding with the error provided by the standard \gls{rom} and a full Gauss quadrature rule determined by the \gls{fe} mesh in the microcells.

Note that by taking a small number of strain modes (\num{20} modes) and cubature points (less than \num{400}), the assessed errors for both cases (\SI{0.2}[<]{\percent}) are much lower than the typical admissible tolerances in engineering problems.
\FloatBarrier
%%%%%%%%%%%%%%%%%%%%%%%%%%%%%%%%%%%%%%%%%%%%
\subsection{Model conservation: assessing the physics of the \emph{\gls{hprfe2}} reduced model}
\label{sect:completeness}
%%%%%%%%%%%%%%%%%%%%%%%%%%%%%%%%%%%%%%%%%%%%
%
The \gls{roec} scheme used to integrate the balance equation \eqref{equilibrium_PI_2-1} requires the computation of stresses and the damage evolution \eqref{damage evolution_PI_3-1} in the specific points of the microcell, determined by the \gls{roec} technique.
This feature of the \gls{hprfe2} model confers a physical basis with an inherent capacity for modeling material responses not specifically considered during the sampling stage, but that still obey basic physical assumptions made at the microscale (e.g., constitutive model thermodynamic consistency, loading-unloading conditions).
In this sense, the physics of the reduced model is identical to that shown by the conventional \gls{fe} approach, where the material physical laws are strictly satisfied in the quadrature points.

To test this feature, the microcells of the models \emph{A} and \emph{B} are solved with the \gls{hprfe2} technique using \num{30} strain modes and \num{400} cubature points.
A cyclic macrostrain trajectory is imposed, defined by the macrostrain \eqref{eq:validation_trajectory} and the parameter $\chi$ increasing monotonically from $0$ to $\chi_\text{end}$ and back to $0$.

\begin{figure}[hbt]
    \centering
    \begin{subfigure}[b]{0.4\textwidth}
        \centering
        \includegraphics[width=0.99\textwidth]{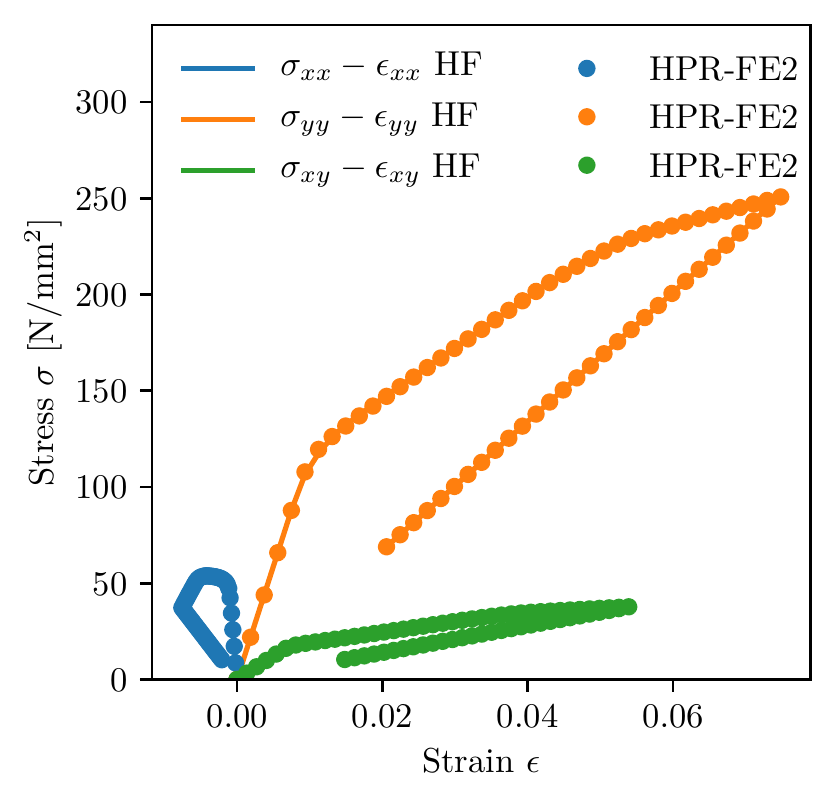}
        \caption{}
        \label{fig_reconstruction_ref_a}

    \end{subfigure}
    \centering
    \begin{subfigure}[b]{0.59\textwidth}
       \centering{
        \begin{subfigure}[b]{0.6\textwidth}
        \centering
            \includegraphics[height=3.0cm]{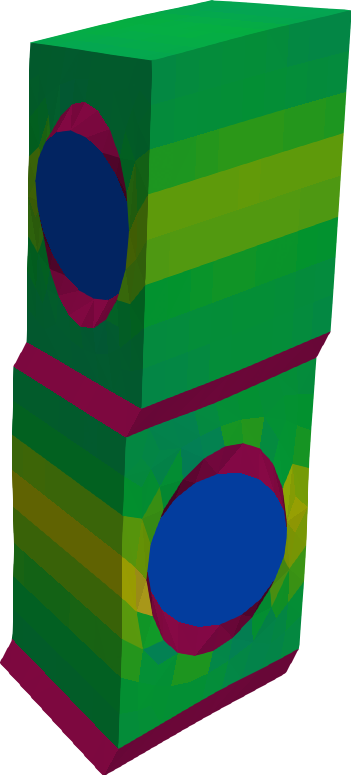}
            \includegraphics[height=3.0cm]{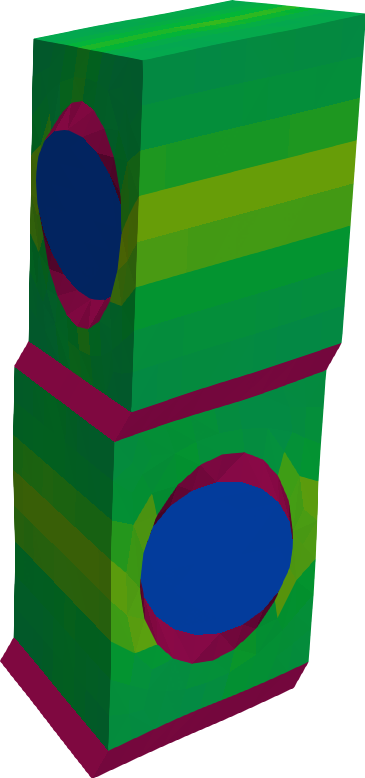}
            \includegraphics[height=3.0cm]{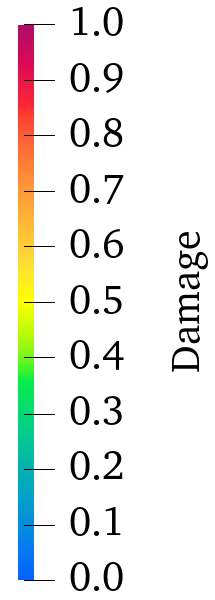}

            \includegraphics[height=3.0cm]{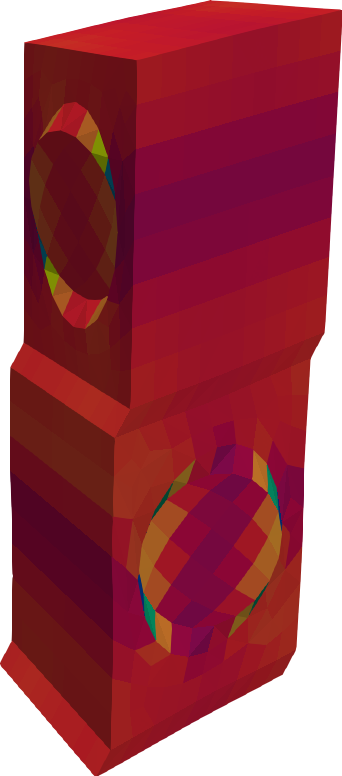}
            \includegraphics[height=3.0cm]{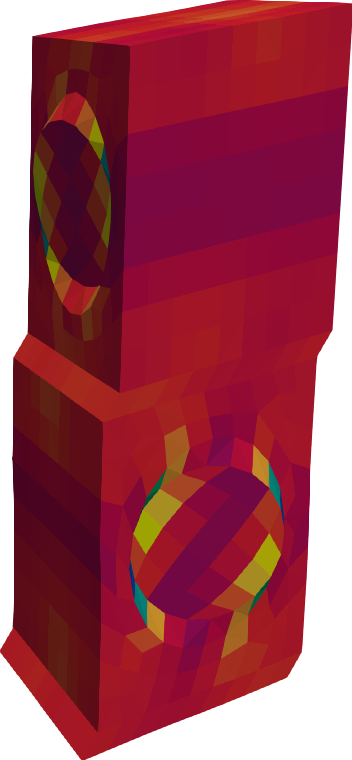}
            \includegraphics[height=3.0cm]{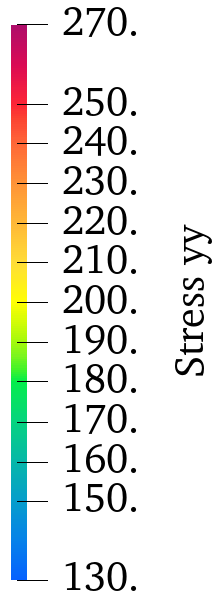}
        \end{subfigure}
        \begin{subfigure}[b]{0.10\textwidth}
            \includegraphics[width=0.99\textwidth]{ref_local_4}
        \end{subfigure}
        \caption{}
        \label{fig_reconstruction_ref_b}
        }
    \end{subfigure}
    \caption{{\bf Model \emph{A}} a) Effective (homogenized) stress vs. effective (macro) strain. Comparison of \gls{hf} and \gls{hprfe2} model solutions in a non-sampled cyclic trajectory.
		b) Top: \gls{hf} (left) and reconstructed (right) damage fields (maximum reconstruction error: \SI{3.24}{\percent}).
        Bottom: \gls{hf} (left) and reconstructed (right) microstress $\sigma_{{\mu}, yy}$ (maximum reconstruction error: \SI{2.60}{\percent}).
		Deformed microcells resulting from the total displacement field reconstruction. Fields and errors correspond to maximum load.
		}
    \label{fig_reconstruction_ref}
\end{figure}

\begin{figure}[hbt]
    \centering
    \begin{subfigure}[b]{0.4\textwidth}
        \centering
        \includegraphics[width=0.99\textwidth]{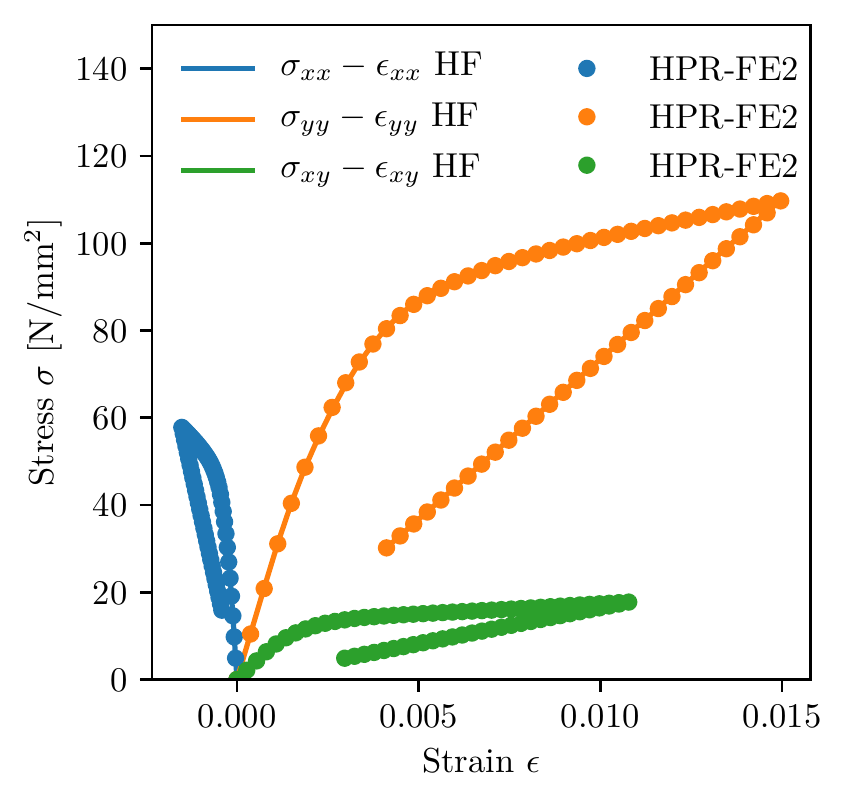}
        \caption{}
        \label{fig_reconstruction_porto_a}
    \end{subfigure}
    \centering
    \begin{subfigure}[b]{0.59\textwidth}
        \centering{
        \begin{subfigure}[b]{0.65\textwidth}
        \centering
        \includegraphics[height=3.0cm]{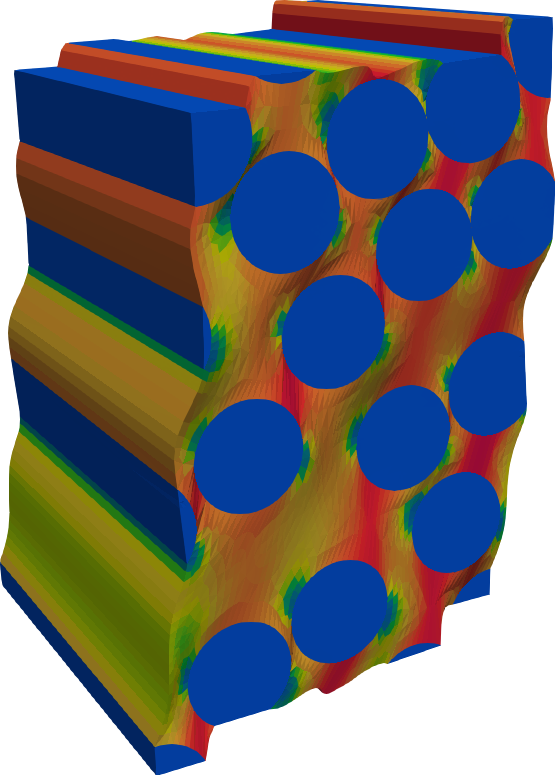}
        \includegraphics[height=3.0cm]{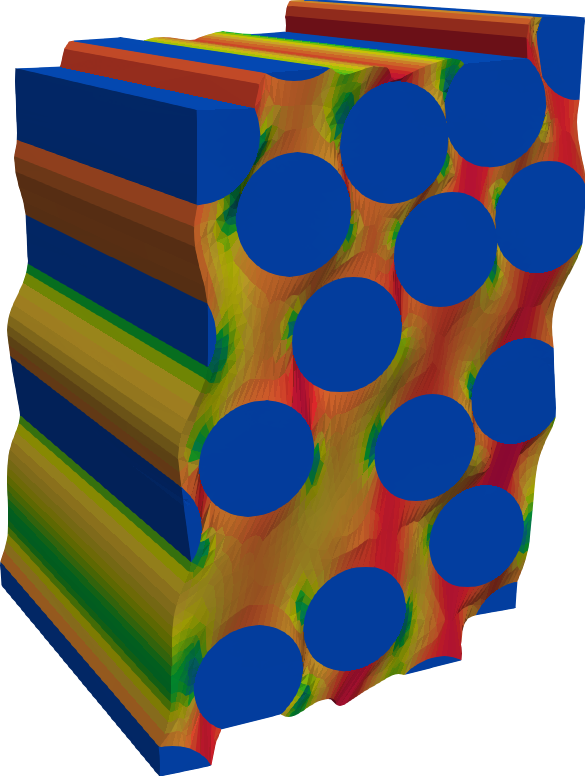}
        \includegraphics[height=3.0cm]{damage_bar_porto}

        \includegraphics[height=3.0cm]{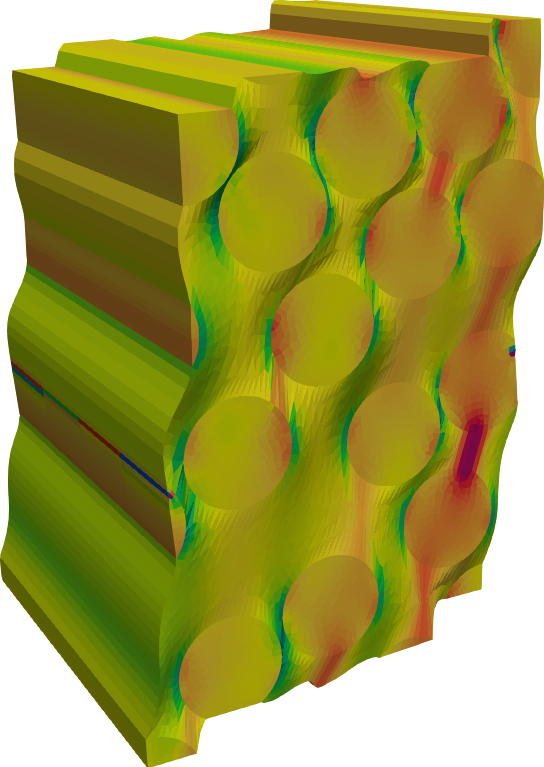}
        \includegraphics[height=3.0cm]{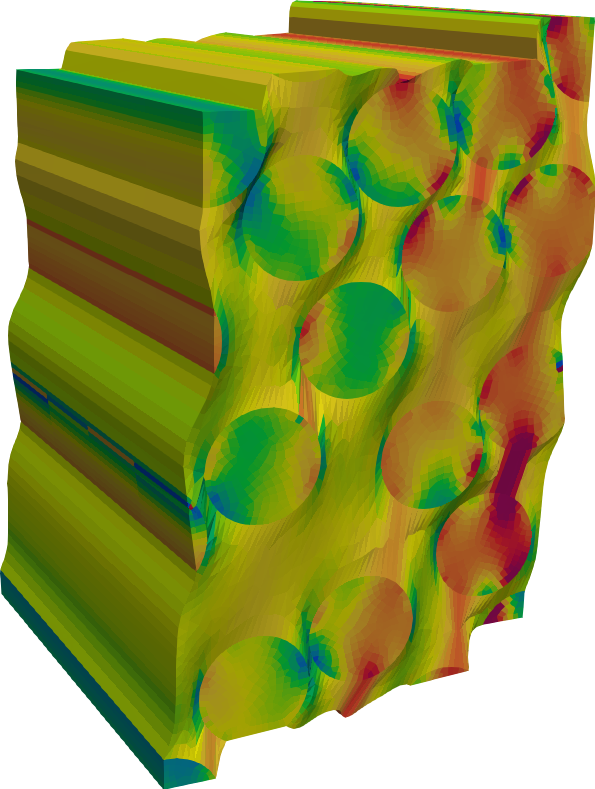}
        \includegraphics[height=3.0cm]{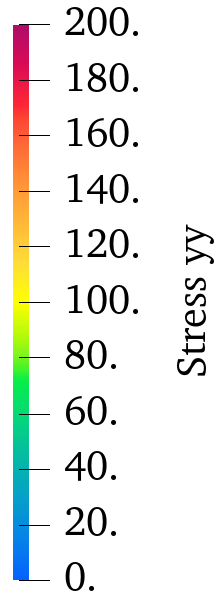}
        \end{subfigure}
        \begin{subfigure}[b]{0.10\textwidth}
            \includegraphics[width=0.99\textwidth]{ref_local_4}

        \end{subfigure}
        \caption{}
        \label{fig_reconstruction_porto_b}
        }
    \end{subfigure}
    \caption{{\bf Model \emph{B}} a) Effective (homogenized) stress vs effective (macro) strain. Comparison of \gls{hf} and \gls{hprfe2} solutions in a non-sampled cyclic trajectory.
		b)  Top: \gls{hf} (left) and reconstructed (right) damage fields (maximum reconstruction error: \SI{6.20}{\percent}).
        Bottom: \gls{hf} (left) and reconstructed (right) microstress $\sigma_{{\mu}, yy}$ (maximum reconstruction error: \SI{6.24}{\percent}).
        Deformed microcells resulting from the total displacement field reconstruction. Fields and errors correspond to maximum strain.
    }
    \label{fig_reconstruction_porto}
\end{figure}

Different homogenized stress components are depicted in Figures~\ref{fig_reconstruction_ref_a} and~\ref{fig_reconstruction_porto_a} for both microcells.
Continuous curves depict solution paths of the \gls{hf} model while dashed curves correspond to the reduced model.
In both cases the solution of the reduced model is remarkably coincident with the \gls{hf} solution even during the unloading regime, which was not tested in any of the sampling paths commented in the previous section.
This result evidences the close connection between the reduced model and the expected physical response of the composite.

Deformed microcells at maximum load $\chi_\text{end}$ are depicted in Figures~\ref{fig_reconstruction_ref_b} and~\ref{fig_reconstruction_porto_b}.
These deformed configurations have been obtained with the reconstructed displacements fields, as indicated in Section~\ref{sect:recovery}.
In the same figures, the reconstructed damage $d_{\mu}$ and stress component $\sigma_{\mu,yy}$ are plotted on the deformed cells.
The reconstructed fields are bounded by an error of  3.24\% for model A, and 6.24\% for model B, compared with the fields obtained with the \gls{hf} technique.
Note the very high accuracy of the \gls{hprfe2} solutions to display the strain jumps caused by the fiber-matrix decohesion and interply delamination mechanisms.

\FloatBarrier
%
%%%%%%%%%%%%%%%%%%%%%%%%%%%%%%%%%%%%%%%%%%%%
\subsection{Computational performance}
\label{sect:speedup}
\begin{figure}
	\centering
	\includegraphics[width=.75\textwidth]{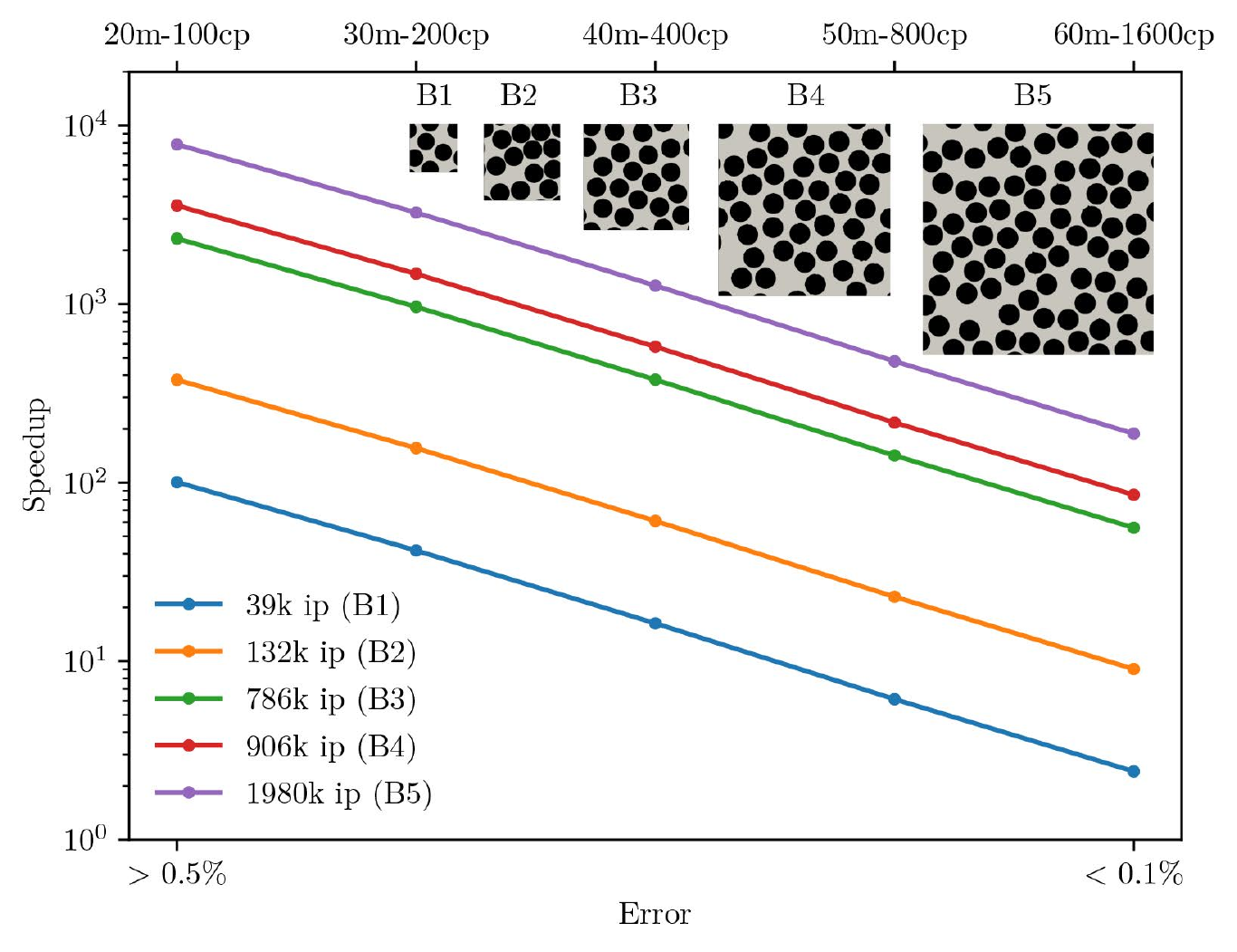}
	\caption{Speedup of microcells B1 to B5 of increasing number of Gauss integration points (ip), for different combinations of number of modes and cubature points. For all the cases, the error range is \SIrange{0.1}{0.5}{\percent}.}
	\label{fig:speedup_modes}
\end{figure}

%%%%%%%%%%%%%%%%%%%%%%%%%%%%%%%%%%%%%%%%%%%%
%
Gains in computing time of \gls{hprfe2} with respect to the standard \gls{fe2} is quantified by the speedup curves shown in Figure~\ref{fig:speedup_modes}.
Speedup curves correspond to 3D microcells B1 to B5, of increasing mesh size
with \numlist{4758; 16672; 99040; 114344; 249880} elements (approximately \SI{96}{\percent} hexahedra and \SI{4}{\percent} 6-ip wedges).
The microcells are based on the described model \emph{B}, and present a random distribution of fibers (volume fraction of \SI{60}{\percent}) with a periodic boundary, generated with the algorithm described in \cite{melro2008generation}.
Each speedup value is computed as the ratio between the time of the \gls{hf} analysis and the corresponding reduced model analysis, varying in the number of modes and cubature points.
\gls{hf} CPU times are \SIlist{171; 640; 3960; 6060; 13350}{s} for models R1 to R5, and CPU times demanded by the reduced model are \SIlist{1.7; 4.1; 10.5; 27.9; 70.9}{s} for combinations 20m-100ip, 30m-200ip, 40m-400ip, 50m-800ip and 60m-1600ip, respectively.
In the case of reduced analysis, time depends only on the number of modes and cubature points, indistinctly of the RVE that they represent.

As expected, the smaller the number of modes and cubature points, the lower the accuracy and the higher the speedup.
Also, the speedup is higher for larger microcells and a given number of strain modes and cubature points.
Consequently, in order to determine a fair measure of the speedup for engineering purposes, the maximum admissible error of the evaluated results has to be fixed as a target, which in turn limits the minimum number of strain modes and cubature points that can be taken in the reduced model.

\begin{figure}
    \centering
    \includegraphics[width=.65\textwidth]{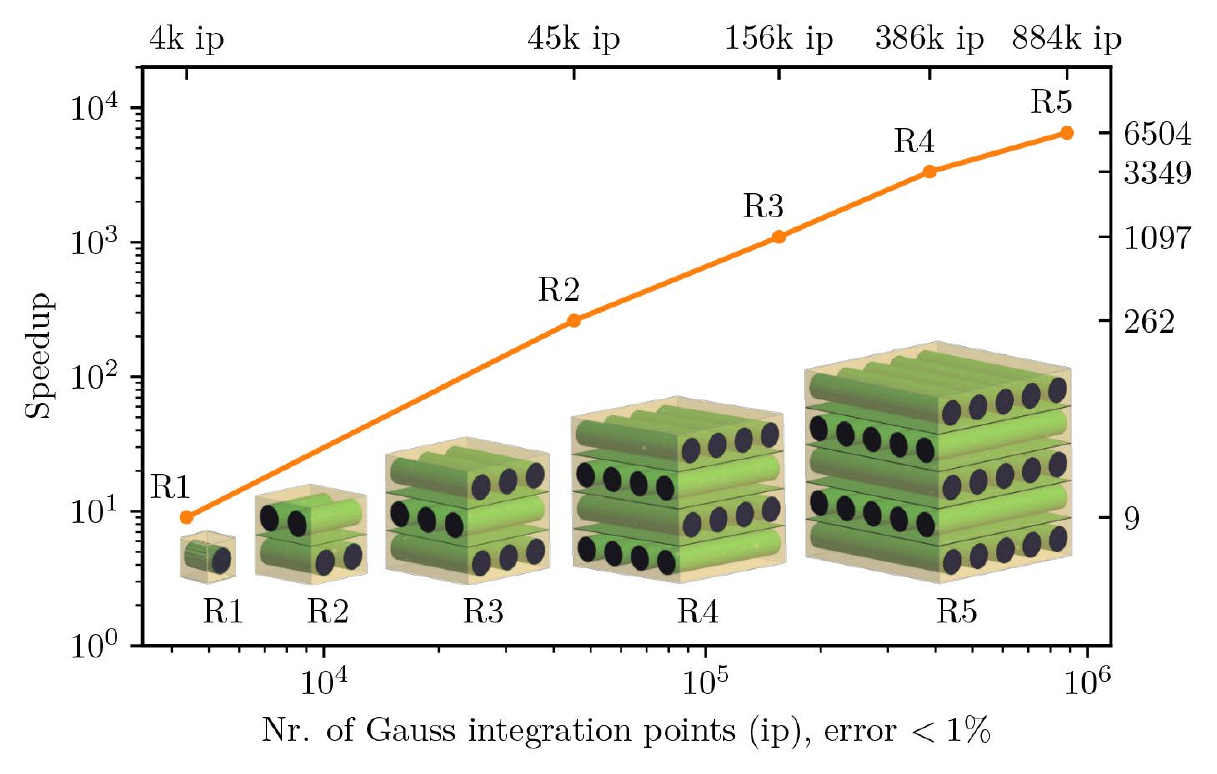}
    \caption{Speedup values obtained for a sequence of microscale problems involving microcells R1 to R5, of increasing complexity (evaluated in terms of the number of Gauss integration points of the \gls{fe} discretization).  The maximum error allowed is \SI{1}{\percent}.}
    \label{fig:speedup}
\end{figure}

In this view, a series of tests involving microcells of increasing
mesh size (\numlist{546; 5656; 9446; 48285; 110480} hexahedra in microcells R1 to R5, respectively)
and complexity (defined as the number of algebraic operations) has been previously performed by the authors in~\cite{Lloberas_Complas_2019}, and the resulting speedups are reproduced in Figure~\ref{fig:speedup}
(involved CPU times are shown later in section~\ref{sect:dogbone_times}).
For every case of the study, the errors obtained by reduced models when taking at least \num{40} strain modes and \num{200} cubature points are below \SI{1}{\percent} (an acceptable error threshold for engineering purposes).
Note that periodicity of the solutions is specifically broken by the fact that the geometry is not repetitive, i.e., the larger microcells are not the tiling of the smaller ones.
Furthermore, boundary conditions compatible with minimal kinematic constraints are defined, so that the physical complexity increases with the volume (e.g., the strain field of the 4-fiber-4-layer microcell R4 is not the periodic repetition of the strain field of the 2-fiber-2-layer microcell R2.)
Thus, the strain modes of the larger microcells are not the periodic repetition the smaller ones.

It is remarkable that, for large microcells with millions of Gauss integration points within the \gls{fe} discretization, the speedup of the reduced model reaches values in the order of \num{e4} (plus any other speedup gain due to parallel processing).
In other words, \gls{hf} multiscale simulations that would take years of computing time, can now be performed in a few hours with errors below \SI{1}{\percent}.

\FloatBarrier
%
%%%%%%%%%%%%%%%%%%%%%%%%%%%%%%%%%%%%%%%%%%%%
\subsection{Material customization: dependency on material properties values for the same material model}
\label{sect:customization}
%%%%%%%%%%%%%%%%%%%%%%%%%%%%%%%%%%%%%%%%%%%%
%
Several tests are outlined to prove the customization of the methodology in terms of changes in the material parameters characterizing the composite phases.
The idea behind these numerical experiments is to evaluate the capability of the \gls{hprfe2} technique---specifically developed using a  given composite morphology and a given set of material parameters---for simulating composites having the same morphology and ruled by the same constitutive model, but with remarkably different values of material properties with respect to those adopted for the sampling stage.
A suitable model response for this kind of test may be a highly desired attribute in applications of material design, as well as in virtual testing.

The microcells employed in the previous sections have been constructed with sampling procedures that have considered the reference material parameters depicted in Tables~\ref{tab:material_params_ref} and~\ref{tab:material_params_porto}.
In the following numerical experiments, the material response of those reference materials are compared to the same reduced models, but in this case using a different set of material properties without performing new sampling.
The microcells of the models \emph{A} and \emph{B} are now characterized by the matrix and cohesive interphase parameters displayed in Table \ref{tab:custom_params} and denoted materials \emph{Custom~M1}, \emph{Custom~M2}.
Note that \emph{Custom~M3} is only used in the case of the microcell of the model \emph{B}, in which the fibers are subjected to damage.
 In the table, the new material parameters are compared to the reference parameters used in the sampling of the reduced models.

\begin{figure}[hbt]
    \centering
    \begin{subfigure}[hbt]{0.49\textwidth}
        \includegraphics[width=\textwidth]{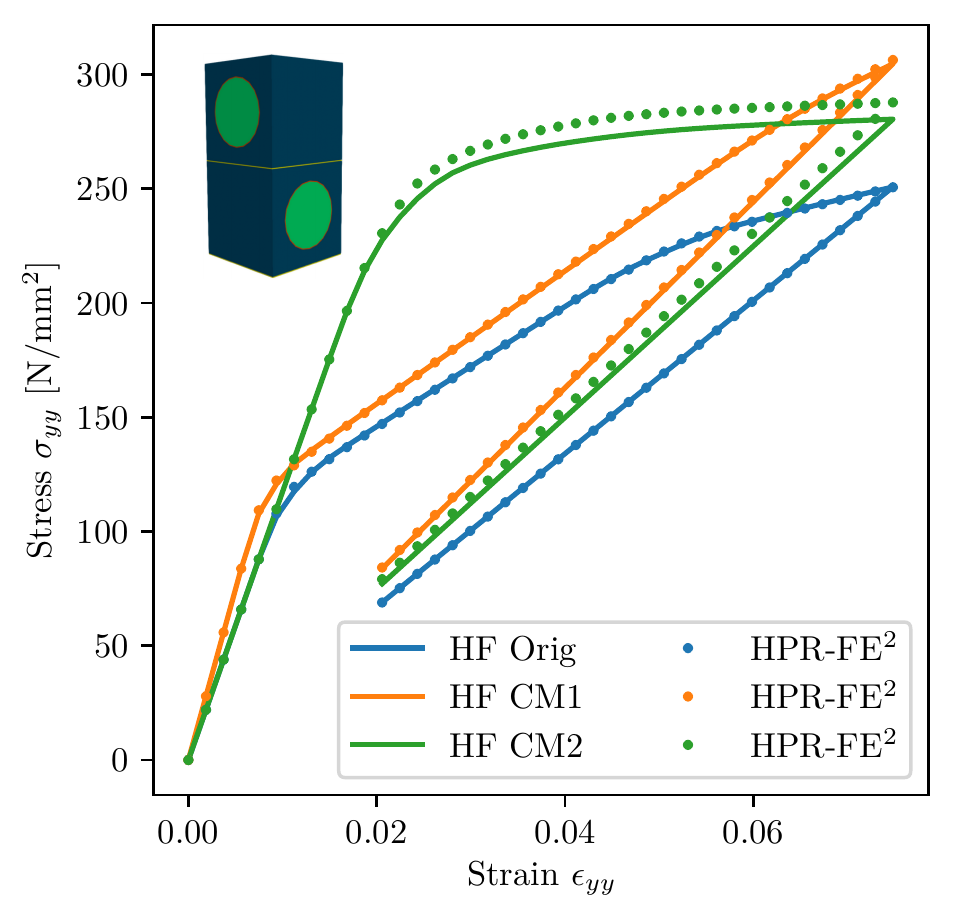}
    \caption{Model \emph{A}}
    \label{fig:custom_params_a}
    \end{subfigure}
    \begin{subfigure}[hbt]{0.49\textwidth}
        \includegraphics[width=\textwidth]{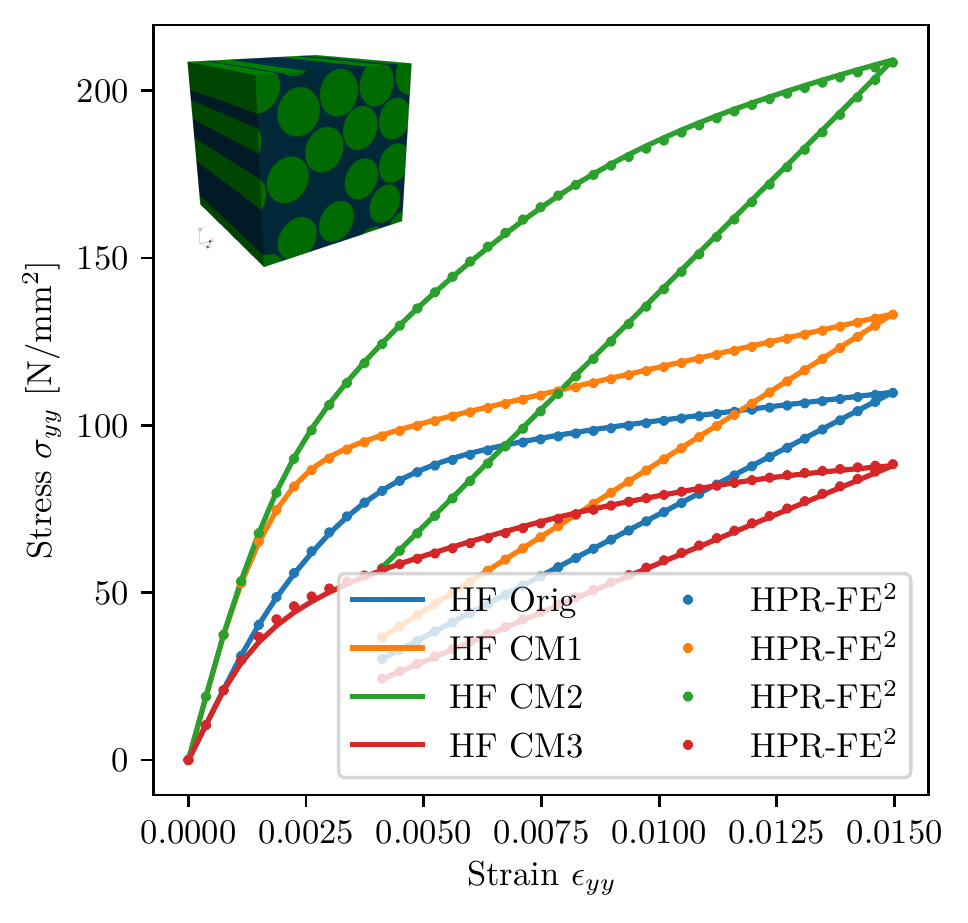}
    \caption{Model \emph{B}}
    \label{fig:custom_params_b}
    \end{subfigure}
    \caption{%
        Numerical experiments with customized materials.
        Effective stress $\sigma_{yy}$ versus effective strain $\varepsilon_{yy}$.
        Solutions of customized tests are displayed for the \gls{hf} and the \gls{hprfe2}.
        Original material parameters used for the \gls{hprfe2} model sampling; customized material parameters as shown in Table~\ref{tab:custom_params}. }
    \label{fig:custom_params}
\end{figure}

Custom reduced model responses are compared to those obtained using \gls{hf} models.
Result curves $\sigma_{yy}$ versus $\varepsilon_{yy}$ for each custom material are reported in Figure~\ref{fig:custom_params}, for \gls{hprfe2} bases with \num{30} strain modes and \num{400} cubature points.
Curves correspond to an imposed cyclic macrostrain trajectory, defined by the macrostrain \eqref{eq:validation_trajectory}, with $\chi$ increasing monotonically from $0$ to $\chi_\text{end} = 0.1$ (microcell \emph{A}) and  $\chi_\text{end} = 0.02$ (microcell \emph{B}), and back to 0.
The results obtained with the original material parameters used for the sampling, and for a similar trajectory, are also included in the figure.

These plots show a close agreement between the solutions obtained with both methodologies.
A slight difference is observed in Figure~\ref{fig:custom_params_a} for microcell of the model \emph{A}, \emph{Custom M2} material.
Possibly, this effect is due to the extremely large change in the customized hardening law.

Note that in the case of the  curve CM3 in Figure~\ref{fig:custom_params_b}, the fiber is allowed to damage even though the \gls{hprfe2} strain bases used in this analysis has been obtained considering the fibers elastic.
Even so, the comparison with the \gls{hf} is remarkable good.
\begin{table}[hbt]
    \centering
    \begin{tabular}{llc|ccc}
        \toprule
        &       & Original                                              &  Custom M1  &           Custom M2 &           Custom M3                                           \\
        \midrule
        \textbf{Model \emph{A}}&&&&& \\
        Matrix: & Young modulus E $[N/mm^2]$                           & \num{3790}  &         \textbf{5000}         &           \num{3790}  &     --    \\
        & Poisson ratio $\nu$                                  & \num{0.37}  &          \num{0.37}           &           \num{0.37}   & --       \\
        & Elastic threshold $\sigma_{\text{e}}$ $[N/mm^2]$     &  \num{200}  & \textbf{300} &           \num{200} &  --         \\
        & Inf elastic threshold $\sigma_{\infty}$ $[N/mm^2]$ &  \num{220}  & \textbf{320} &           \num{220 }  &  --       \\
        & Hardening parameters $H_1$, $H_2$  & \num{0.10}, \num{0.01}  & \textbf{0.5}, \textbf{0.1} &  \num{0.10},  \num{0.01} & --\\ \midrule
        Cohesive, & Young modulus E $[N/mm^2]$                           & \num{3790}  &          \num{3790}           &          \num{ 3790} &  --        \\
        Interply:& Poisson ratio $\nu$                                  & \num{0.37}  &          \num{0.37}           &          \num{ 0.37}       &  --  \\
        & Elastic threshold $\sigma_{\text{e}}$ $[N/mm^2]$     &  \num{100}  &           \num{100}           &           \num{100}      &  --    \\
        & Inf elastic threshold $\sigma_{\infty}$ $[N/mm^2]$ & \num{100.1} &          \num{100.1}          &          \num{ 100.1}      &  --  \\
        & Hardening parameters $H_1$,  $H_2$    & \num{0.01}, \num{0.01} & \num{0.01}, \num{0.01} & \textbf{0.9}, \textbf{0.9} & --\\
        \midrule
        \textbf{Model \emph{B}} & &      &      &     &    \\
        Matrix:& Young modulus E $[N/mm^2]$                              &  \num{4000}  & \textbf{8000} & \textbf{8000} &  \num{4000}  \\
        & Poisson ratio $\nu$                                     &  \num{0.38}  &  \num{0.38}   &  \num{0.38}   &  \num{0.38}  \\
        & Elastic threshold $\sigma_{\text{e}}$ $[N/mm^2]$        &   \num{60}   &   \num{60}    &   \num{60}    &   \num{60}   \\
        & Inf elastic threshold $\sigma_{\infty}$ $[N/mm^2]$    &   \num{70}   &   \num{70}    & \textbf{140}  &   \num{70}   \\
        & Hardening parameters $H_1$, $H_2$                                & \num{0.335}, \num{0.05}    &  \num{0.335}, \num{0.05}    &  \num{0.335}, \num{0.05}    & \num{0.335}, \num{0.05}  \\ \midrule
        Inclusion:& Young modulus E $[N/mm^2]$                           & \num{231000} & \num{231000}  & \num{231000}  & \num{231000} \\
        & Poisson ratio $\nu$                                  &  \num{0.2}   &   \num{0.2}   &  \num{ 0.2}   &  \num{0.2}   \\
        & Elastic threshold $\sigma_{\text{e}}$ $[N/mm^2]$  & --  &  --  &  --  & \textbf{60}  \\
        & Inf elastic threshold $\sigma_{\infty}$ $[N/mm^2]$ & --  &  --  &  --  & \textbf{70}  \\
        & Hardening parameters $H_1$, $H_2$                            &  --  &  --   & --   &  \textbf{0.01},  \textbf{0.01} \\ \bottomrule
    \end{tabular}
    \caption{Material properties for the customization numerical experiments. Properties of columns \emph{Custom M1}, \emph{Custom M2}, and \emph{Custom M3} are compared to the original properties  (now taken as reference) adopted for the sampling process. Property changes in the customization tests are remarked in boldface.}
    \label{tab:custom_params}
\end{table}

\FloatBarrier

% !TeX root = paper_HRFE2_2020.tex
% !TeX spellcheck = en_US
% !TeX encoding = UTF-8

\FloatBarrier
%
%%%%%%%%%%%%%%%%%%%%%%%%%%%%%%%%%%%%%%%%%%%%
\subsection{Virtual testing of a composite coupon}
\label{sect:dogbone}
\begin{figure}[hbt]
    \begin{subfigure}[b]{0.49\textwidth}
        \centering
        \includegraphics[width=0.75\textwidth]{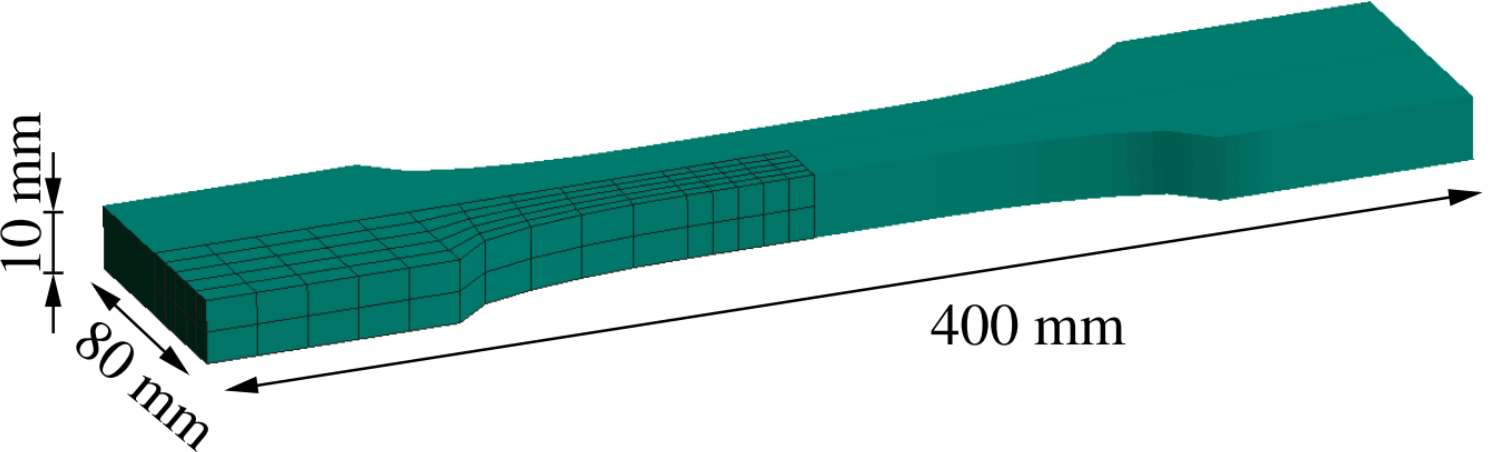}
        \includegraphics[width=0.13\textwidth]{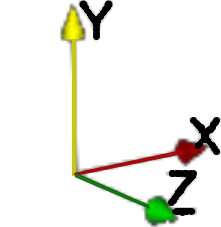}
        \caption{Macroscopic geometry}
        \label{fig:dogbone_model}
    \end{subfigure}
    \begin{subfigure}[b]{0.49\textwidth}
        \centering
        \includegraphics[width=0.6\textwidth]{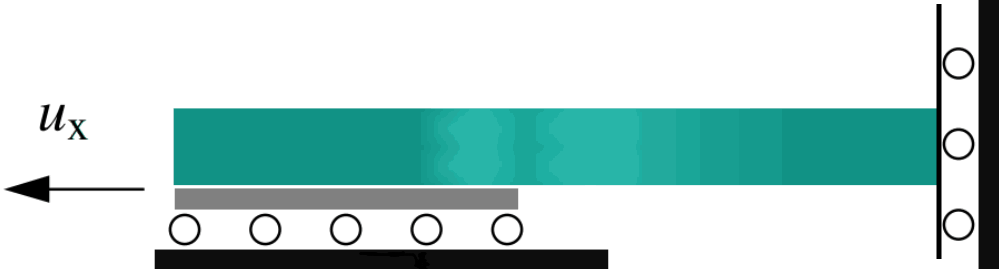}
        \includegraphics[width=0.13\textwidth]{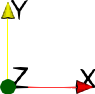}
        \caption{Boundary and load conditions.}
        \label{fig:dogbone_load}
    \end{subfigure}
    \qquad
    \caption{
        Composite coupon under imposed shearing displacement.
        Uniform displacement $u_x$ is transmitted to the coupon using rigid strips acting only on the bottom faces.}
    \label{fig:dogbone}
\end{figure}

%%%%%%%%%%%%%%%%%%%%%%%%%%%%%%%%%%%%%%%%%%%%
%
The coupon depicted in Figure~\ref{fig:dogbone_model} is used to perform several virtual tests.
It is composed by a unidirectional composite ply, characterized by the microcell \emph{B}.
Three cases are simulated, corresponding to fiber orientations \ang{0} (longitudinal), \ang{45} and \ang{90} (transversal) in the plane $x-z$,  obtained by standard rotation of  the microcell to adjust local and global reference axis.
The thickness of the specimen is three orders of magnitude larger  than the fiber diameter, so the separation of scales required by \gls{fe2} method is ensured.

A uniform displacement $u_x$ is imposed on a stiff plate stuck to the inferior surface of the specimen, as shown in Figure \ref{fig:dogbone_load}.
Exploiting the symmetry of the problem, a fourth part of the coupon is simulated with \num{150} linear hexahedral elements and \num{8} Gauss integration points per element.
Each integration point is associated with the microcell described previously in Table \ref{tab:material_params_porto}.

\begin{figure}[hbt]
    \centering
    \begin{subfigure}[b]{0.49\textwidth}
        \includegraphics[width=\textwidth]{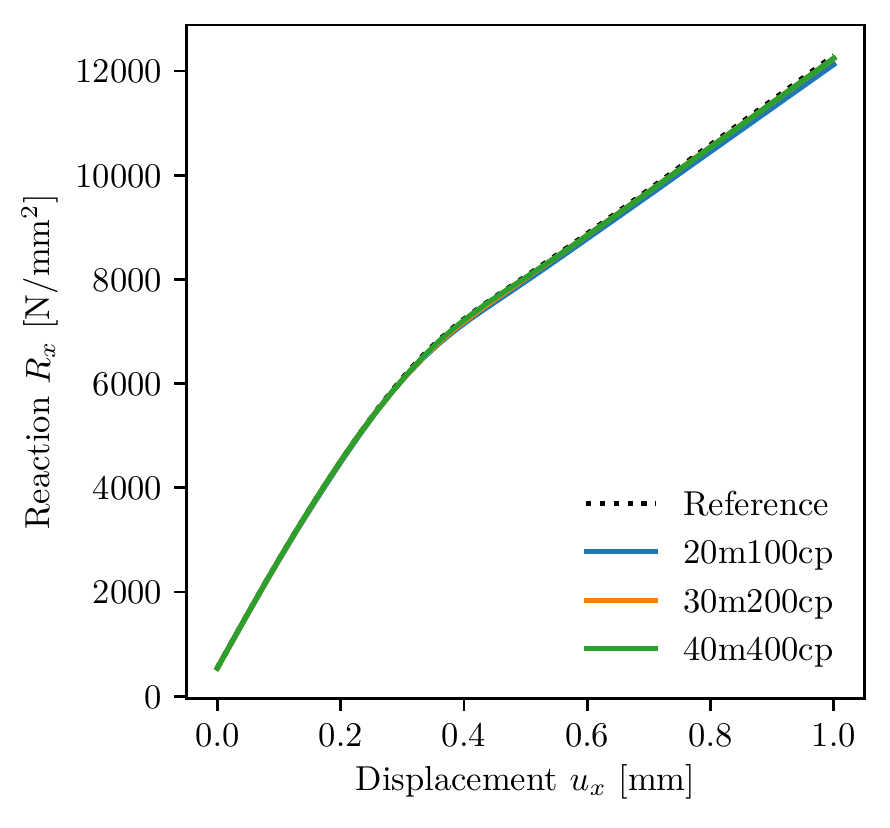}
        \caption{\ang{0}}
    \end{subfigure}
    \begin{subfigure}[b]{0.49\textwidth}
        \includegraphics[width=\textwidth]{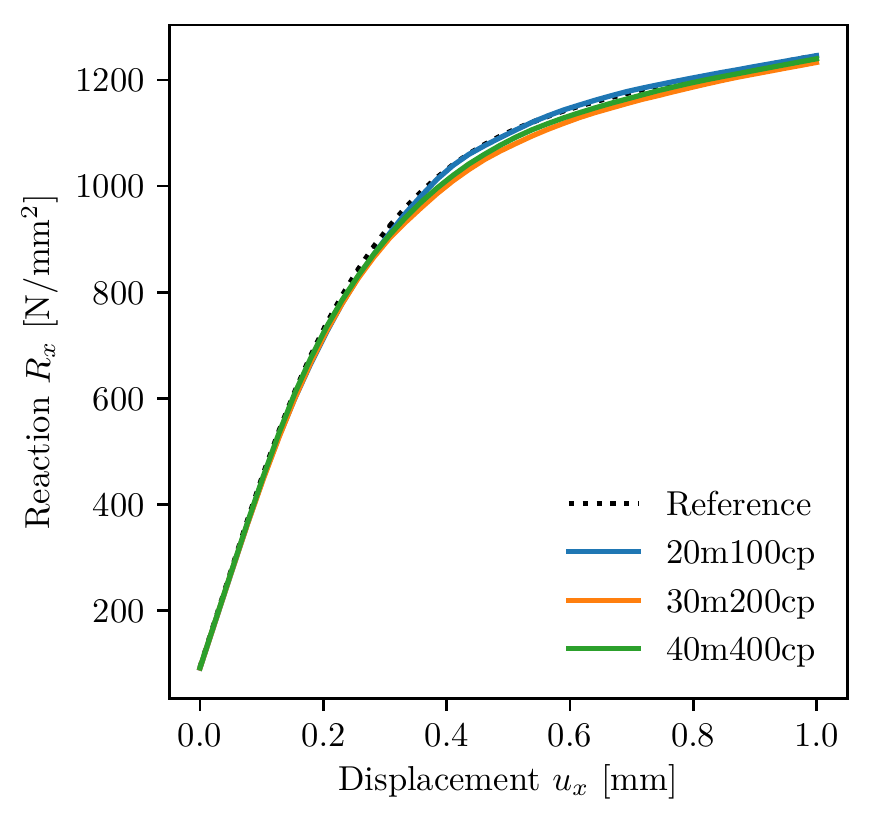}
        \caption{\ang{45}}
    \end{subfigure}
    \begin{subfigure}[b]{0.49\textwidth}
        \includegraphics[width=\textwidth]{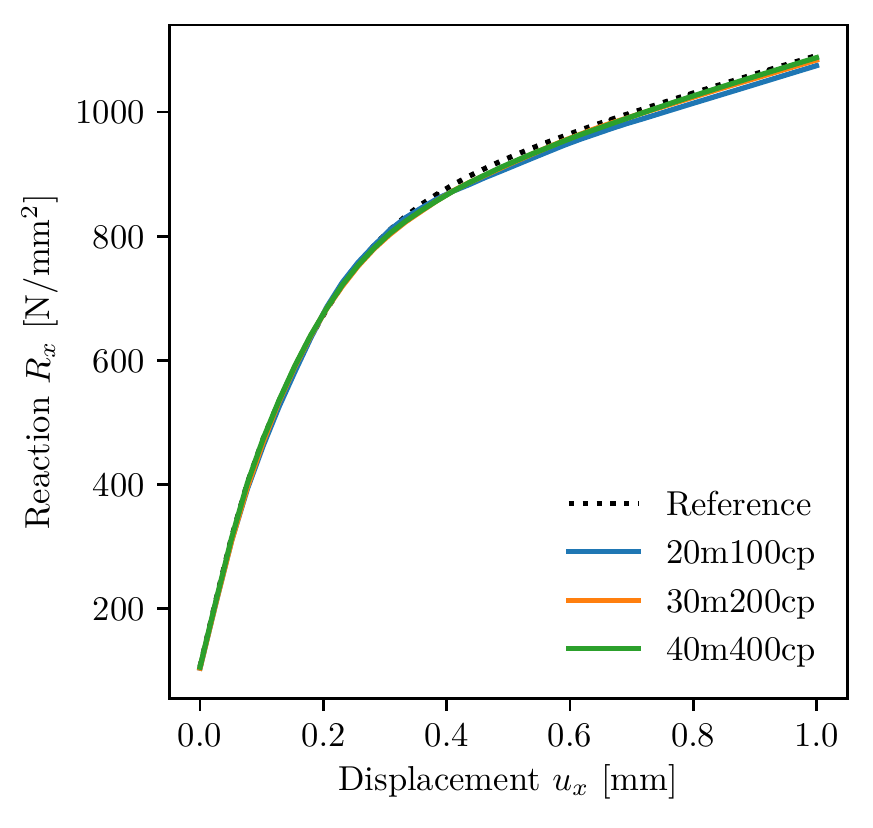}
        \caption{\ang{90}}
    \end{subfigure}
    \caption{
        Response of the composite coupon under axial tension with shear. Reaction force $R_x$ versus displacement $u_x$ for several combinations of modes and cubature points, with fibers orientations \ang{0}, \ang{45} and \ang{90}.
        %
        %\nota{Errors:
         %   (a) \SI{1.1}{\percent}, \SI{0.4}{\percent}, \SI{0.3}{\percent},
         %   (b) \SI{0.5}{\percent}, \SI{1.9}{\percent}, \SI{1.0}{\percent},
         %   (c) \SI{1.4}{\percent}, \SI{0.8}{\percent}, \SI{0.6}{\percent}}
        %
}
\label{fig:dogbone_rve_reaction}
    \end{figure}

The mechanical response of the specimen for different orientations of the fibers is shown in Figure~\ref{fig:dogbone_rve_reaction}.
Solutions of reaction forces versus displacement $u_x$ are plotted for three combinations of microstrain modes and cubature points: 20 modes-100 cubature point (20m100cp), 30 modes-200 cubature points (30m200cp) and 40 modes-400 cubature points (40m400cp).
Due to the unfeasibility of performing an actual \gls{fe2} simulation, the solution taken as reference is the one obtained with the \gls{hprfe2} technique, adopting \num{50} microstrain modes and \num{1800} cubature points.
For the three orientations, the responses of the reduced model agree very well with the reference solutions, with relative errors (with respect to the case 50m-1800cp) between $0.3\%$ and $1.9\%$.
The larger relative error is observed for the case \ang{45}, 30m200cp.
Note that the error, with respect to the reference curve, observed in the case of \ang{45} does not monotonically decrease with the increment of the number of modes and cubature points.
According to the conclusions drawn from Figure~\ref{fig_errors_porto}, we argue that the responses of Figure~\ref{fig:dogbone_rve_reaction} are obtained with a low number of modes and cubature points, in correspondence to a region of non-monotonic error convergence.
However, being that the involved relative errors are small for the purposes of this work, it does not deserve a further analysis involving increments of the number of modes and cubature points.

\begin{figure}
    \begin{subfigure}[b]{\textwidth}
        \begin{subfigure}[b]{0.3\textwidth}
            \centering
            \includegraphics[width=.99\textwidth]{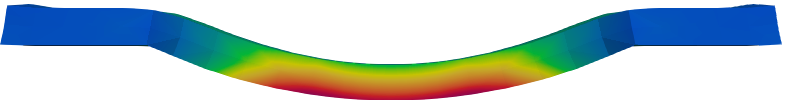}
            \includegraphics[width=.99\textwidth]{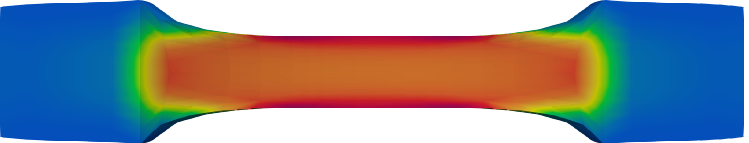}
            \includegraphics[width=.95\textwidth]{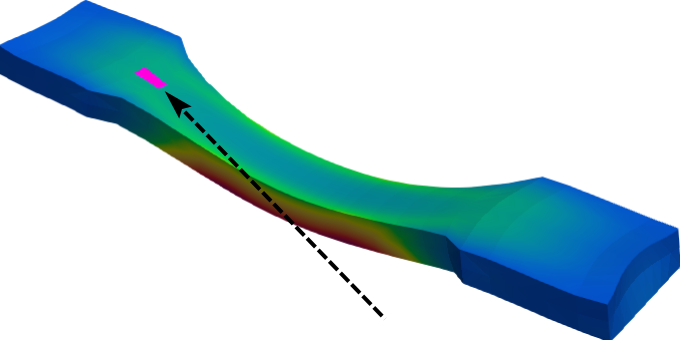}
            \includegraphics[width=0.2\textwidth]{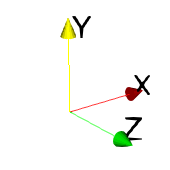}
            \includegraphics[width=0.5\textwidth]{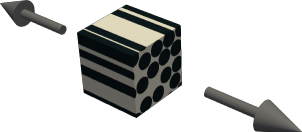}
            \includegraphics[width=\textwidth]{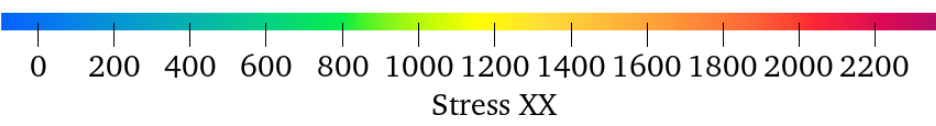}
        \end{subfigure}
        \begin{subfigure}[b]{0.3\textwidth}
            \centering
            \includegraphics[width=.99\textwidth]{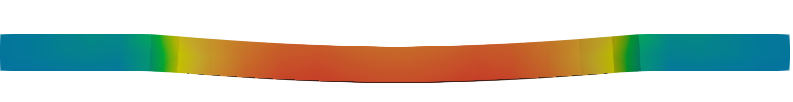}
            \includegraphics[width=.99\textwidth]{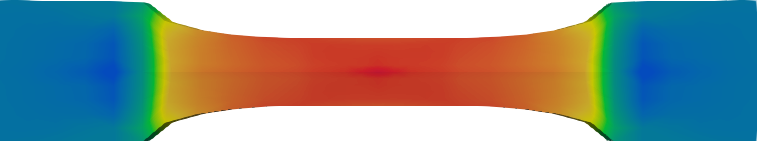}
            \includegraphics[width=.95\textwidth]{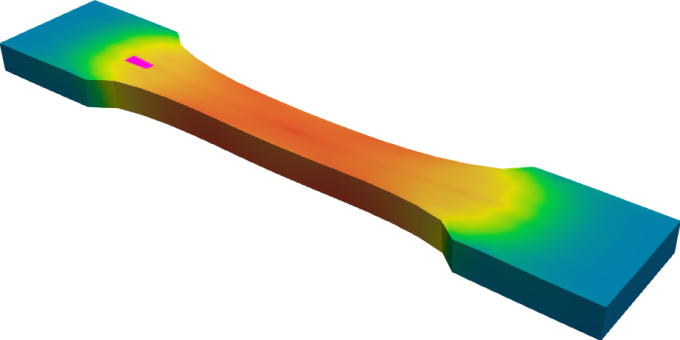}
            \includegraphics[width=0.2\textwidth]{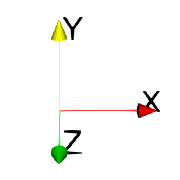}
            \includegraphics[width=0.46\textwidth]{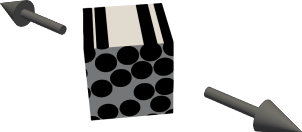}
            \includegraphics[width=0.99\textwidth]{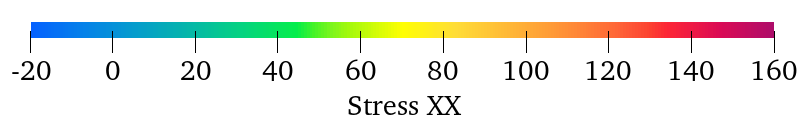}
        \end{subfigure}
        \begin{subfigure}[b]{0.3\textwidth}
            \centering
            \includegraphics[width=0.99\textwidth]{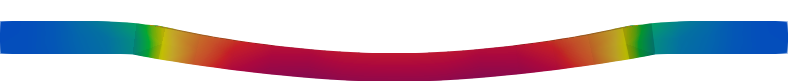}
            \includegraphics[width=0.99\textwidth]{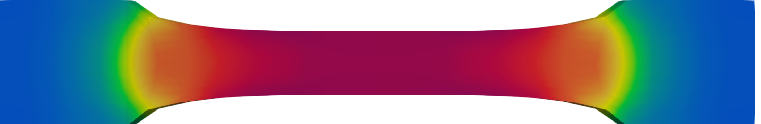}
            \includegraphics[width=.95\textwidth]{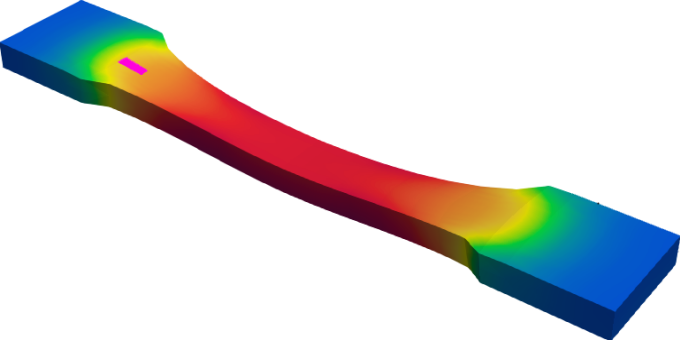}
            \includegraphics[width=0.2\textwidth]{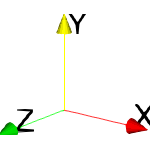}
            \includegraphics[width=0.5\textwidth]{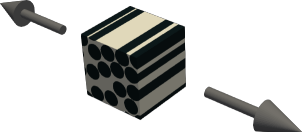}
            \includegraphics[width=\textwidth]{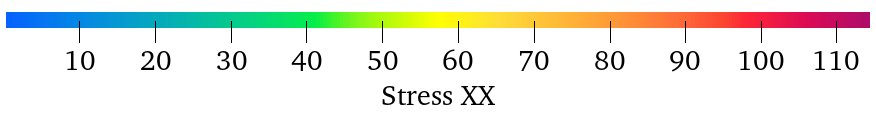}
        \end{subfigure}
        \begin{subfigure}[b]{0.08\textwidth}
    \centering
    \includegraphics[width=0.5\textwidth]{ref_1}
    \includegraphics[width=0.5\textwidth]{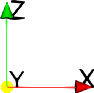}
    \includegraphics[width=0.9\textwidth]{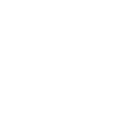}
    \includegraphics[width=0.8\textwidth]{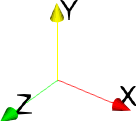}
    \includegraphics[width=0.9\textwidth]{blank}
    \includegraphics[width=0.9\textwidth]{blank}
\end{subfigure}
    \end{subfigure}
    \begin{subfigure}[b]{\textwidth}
        \includegraphics[width=0.25\textwidth]{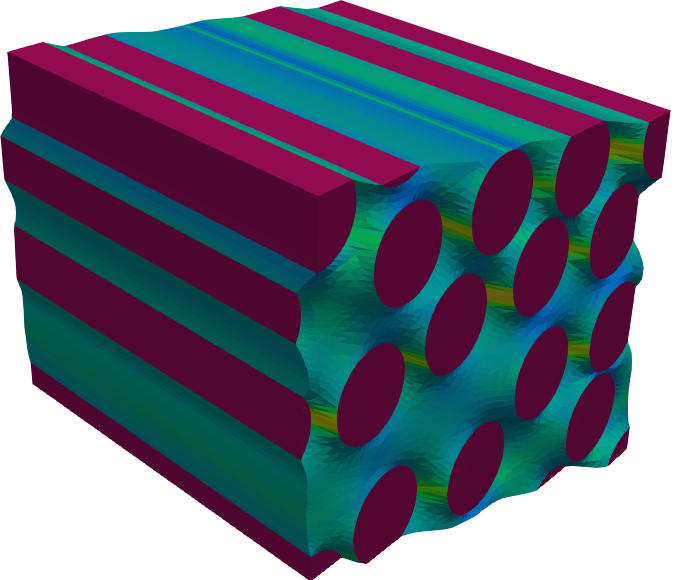}
        \includegraphics[width=0.06\textwidth]{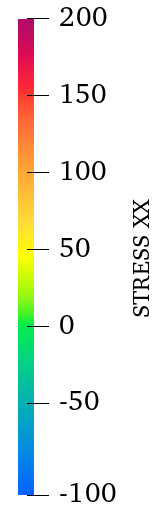}
        \includegraphics[width=0.21\textwidth]{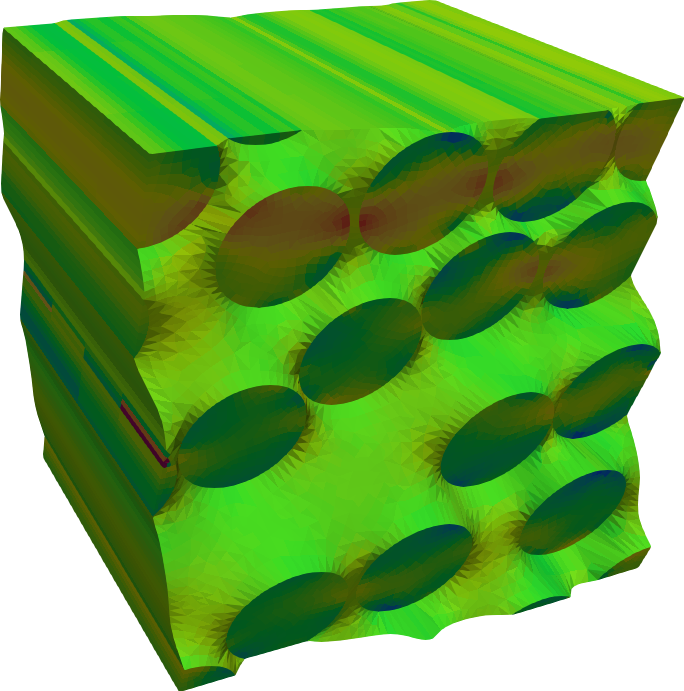}
        \includegraphics[width=0.06\textwidth]{dogbone_rve_stress_bar}
        \includegraphics[width=0.25\textwidth]{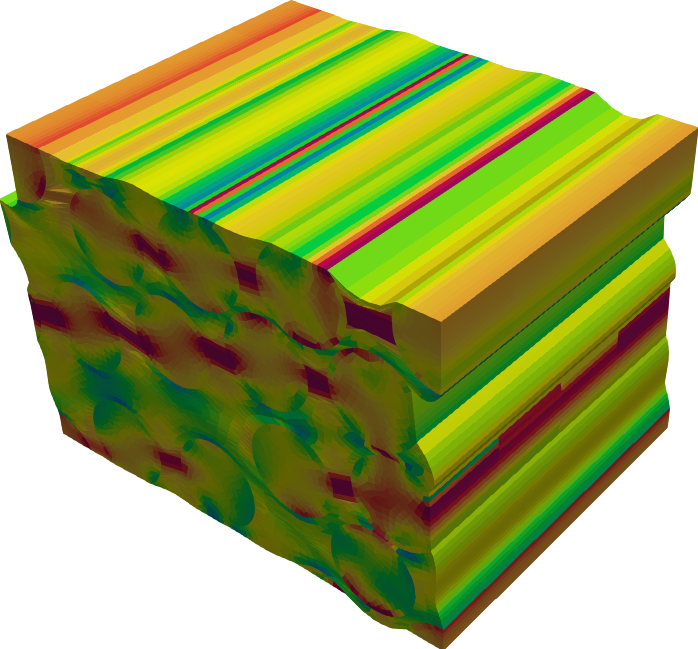}
        \includegraphics[width=0.06\textwidth]{dogbone_rve_stress_bar}
    \end{subfigure}
    \begin{subfigure}[b]{\textwidth}
        \centering
        \includegraphics[width=0.25\textwidth]{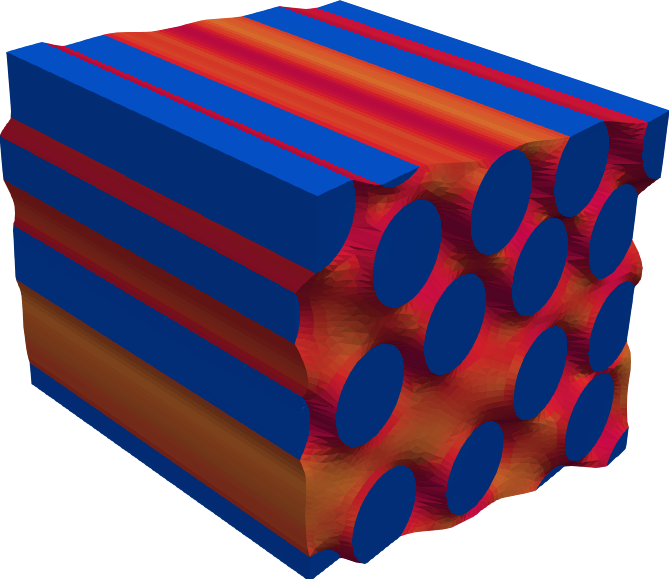}
        \includegraphics[width=0.06\textwidth]{damage_bar_porto}
        \includegraphics[width=0.21\textwidth]{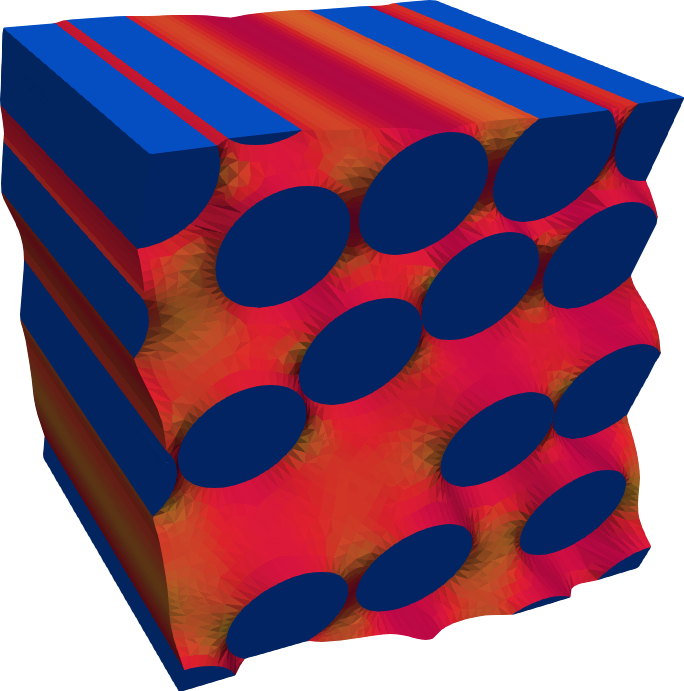}
        \includegraphics[width=0.06\textwidth]{damage_bar_porto}
        \includegraphics[width=0.25\textwidth]{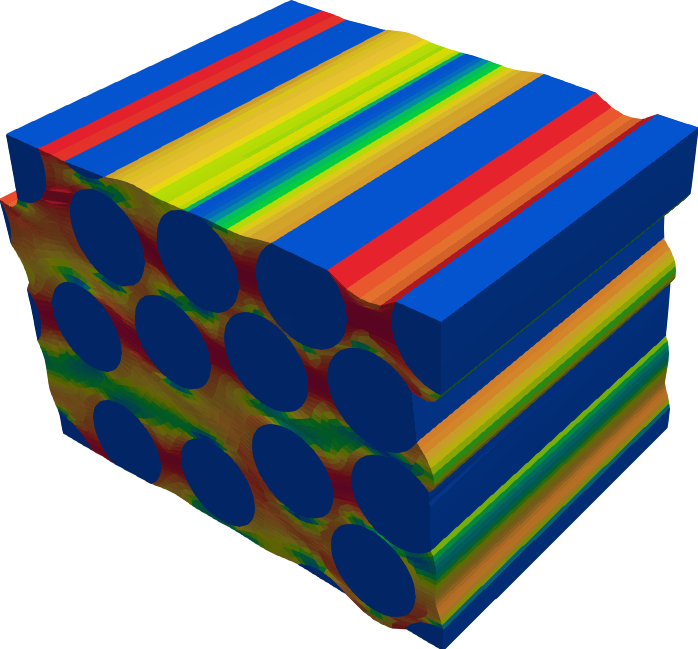}
        \includegraphics[width=0.06\textwidth]{damage_bar_porto}
    \end{subfigure}
    \caption{
        Mechanical response of the composite coupon under axial tension with shear,  for fiber orientations \ang{0}, \ang{45}, and \ang{90} (left, center and right columns, respectively.), for maximum load. Reconstructed fields corresponds to the microcell located at the highlighted element pointed by the arrow.
        On each column, and from top to bottom:
        1) stress field $\sigma_{X\!X}$ of the deformed coupon;
        2) reconstructed fields of the microcell: microstress $\sigma_{\mu, \, xx}$ and damage.
        The fields are plotted onto the deformed microcell resulting from the reconstructed displacement fluctuation field.}
    \label{fig:dogbone_rve}
\end{figure}

Figure~\ref{fig:dogbone_rve} shows the displacement and stress $\sigma_{X\!X}$ of the coupon in the mentioned cases at maximum load.
It also shows the reconstructed microstress $\sigma_{\mu,\,xx}$ and damage field $d_{\mu}$ of the deformed microcell (with the procedure explained in Section~\ref{sect:recovery}), located at an arbitrary Gauss integration point of the coupon.
These microscopic results make evident one of the main features of the \gls{fe2} approaches, and preserved in the \gls{hprfe2} model, i.e., the detailed analysis the microstructural field, which can not  be assessed by simpler phenomenological macroscopic models.
\FloatBarrier
%
%%%%%%%%%%%%%%%%%%%%%%%%%%%%%%%%%%%%%%%%%%%%
\subsubsection{Virtual testing of a composite coupon: customized materials}
\label{sect:dogbone_custom}
%%%%%%%%%%%%%%%%%%%%%%%%%%%%%%%%%%%%%%%%%%%%
%
The previous virtual test is performed with customized material parameters, different from the values used during the sampling stage.
Customized parameters correspond to material \emph{Custom~M2} of Table~\ref{tab:custom_params}.

\begin{figure}[htb]
        \centering
    \begin{subfigure}[b]{0.49\textwidth}
        \includegraphics[width=\textwidth]{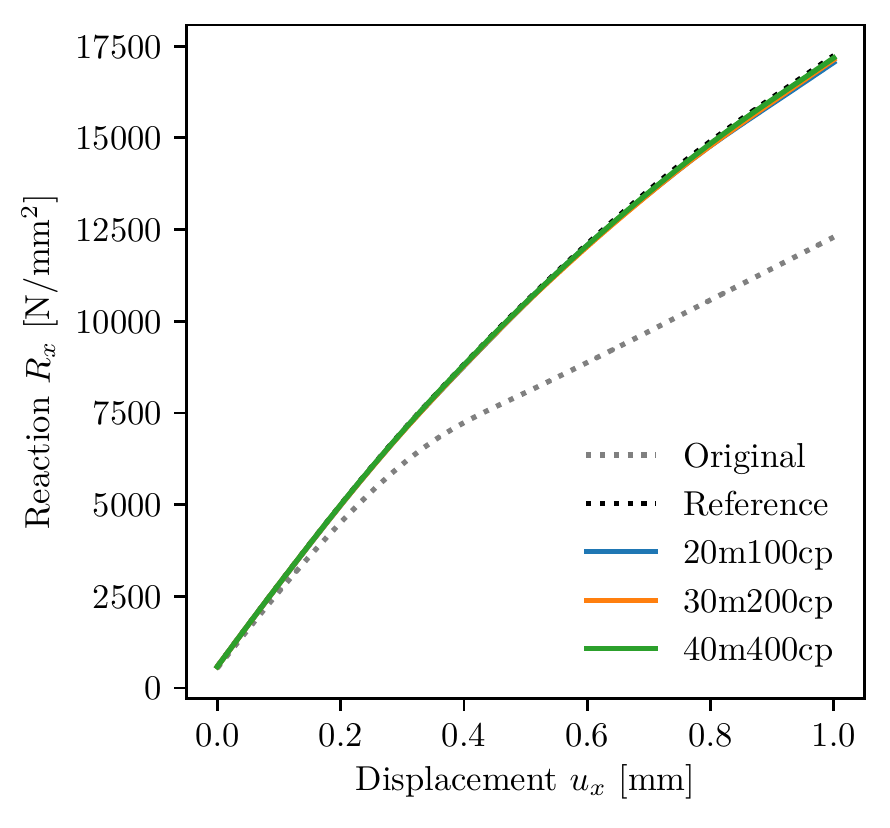}
        \caption{\ang{0}}
    \end{subfigure}
    \begin{subfigure}[b]{0.49\textwidth}
        \includegraphics[width=\textwidth]{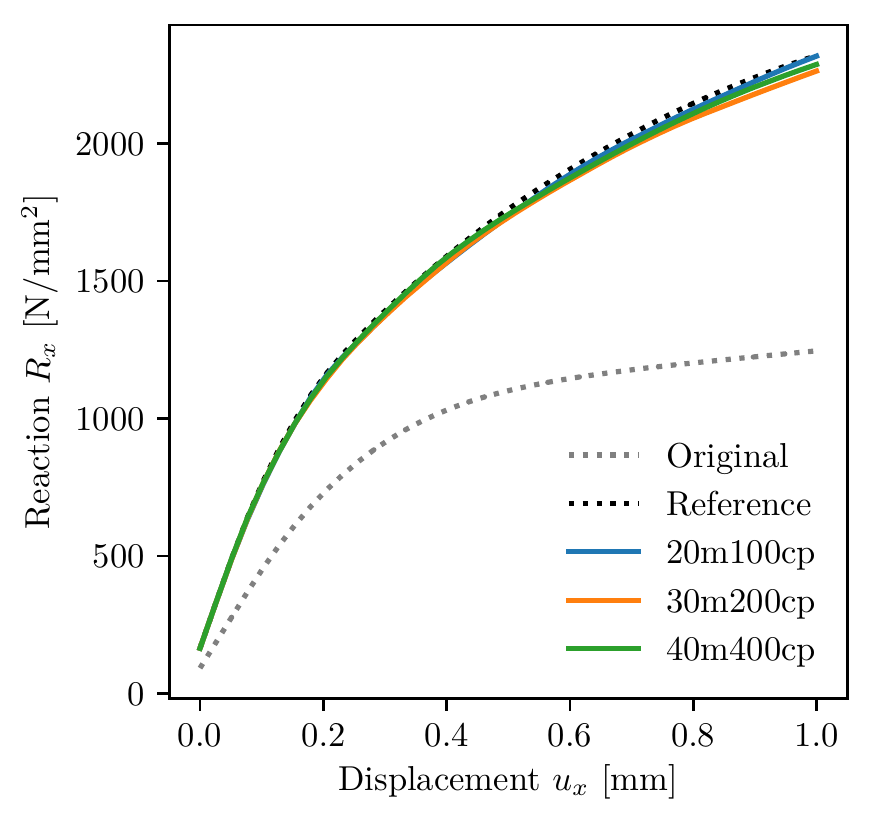}
        \caption{\ang{45}}
    \end{subfigure}
    \begin{subfigure}[b]{0.49\textwidth}
     \includegraphics[width=\textwidth]{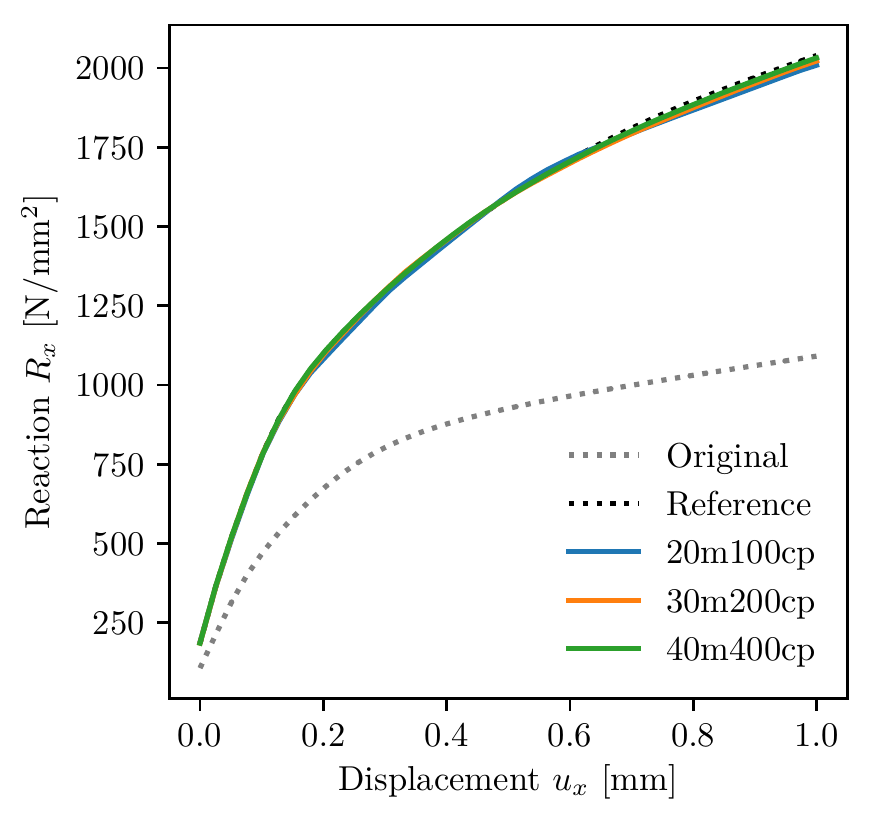}
    \caption{\ang{90}}
    \end{subfigure}
    \caption{Response of a composite coupon under axial tension with shear, with customized material parameters.
        Curves correspond to fiber orientations  \ang{0}, \ang{45} and \ang{90}, and different combinations of number of microstrain modes and cubature points.
        The mechanical response using modified material parameters approximates the ``Reference" curves, obtained from bases  computed with custom materials.
        ``Original" curves are the response curves with the original material parameters.}
    \label{fig:dogbone_reaction_custom}
\end{figure}

The computed results are displayed in the Figure~\ref{fig:dogbone_reaction_custom} for fiber orientations \ang{0}, \ang{45} and \ang{90}.
They  are obtained from the \gls{hprfe2} model sampled with the original material parameters, but using the custom material parameters during the virtual testing stage.
These results are plotted for the cases 20m100cp, 30m100cp and 40m400cp.

The curves denoted as ``Reference" are the result obtained with a \gls{hprfe2} model sampled with the \emph{custom} material parameters.
They are taken as the reference curves.

To compare the structural effects induced by the change of the customized parameters with respect to the original ones adopted in this test, the curves denoted as ``Original" are plotted in the same figure.
These curves, which are the reference curves in Figure~\ref{fig:dogbone_rve_reaction}, are obtained with an \gls{hprfe2} model, sampled with the \emph{original} material parameters and using the same material parameters in the virtual testing stage.

For the three orientations, the responses of the reduced model agree very well with the reference solutions, with relative errors between $0.3\%$ and $2.9\%$.
The larger relative error is observed for the case \ang{45}, 30m200cp.
\FloatBarrier
%
%%%%%%%%%%%%%%%%%%%%%%%%%%%%%%%%%%%%%%%%%%%%
\subsubsection{Virtual testing of a composite coupon: estimated performance comparison}
\label{sect:dogbone_times}
%%%%%%%%%%%%%%%%%%%%%%%%%%%%%%%%%%%%%%%%%%%%
%
The same composite coupon shown in Figure~\ref{fig:dogbone} is tested using a sequence of microcells employed in the speedup curves in Figure~\ref{fig:speedup}.
As discussed before in Section~\ref{sect:speedup}, in terms of mechanical deformation modes, size implies a complexity increment due to their geometry and the use of minimal kinematic boundary conditions.

\begin{figure}[hbt]
	\centering
	\includegraphics[width=.75\textwidth]{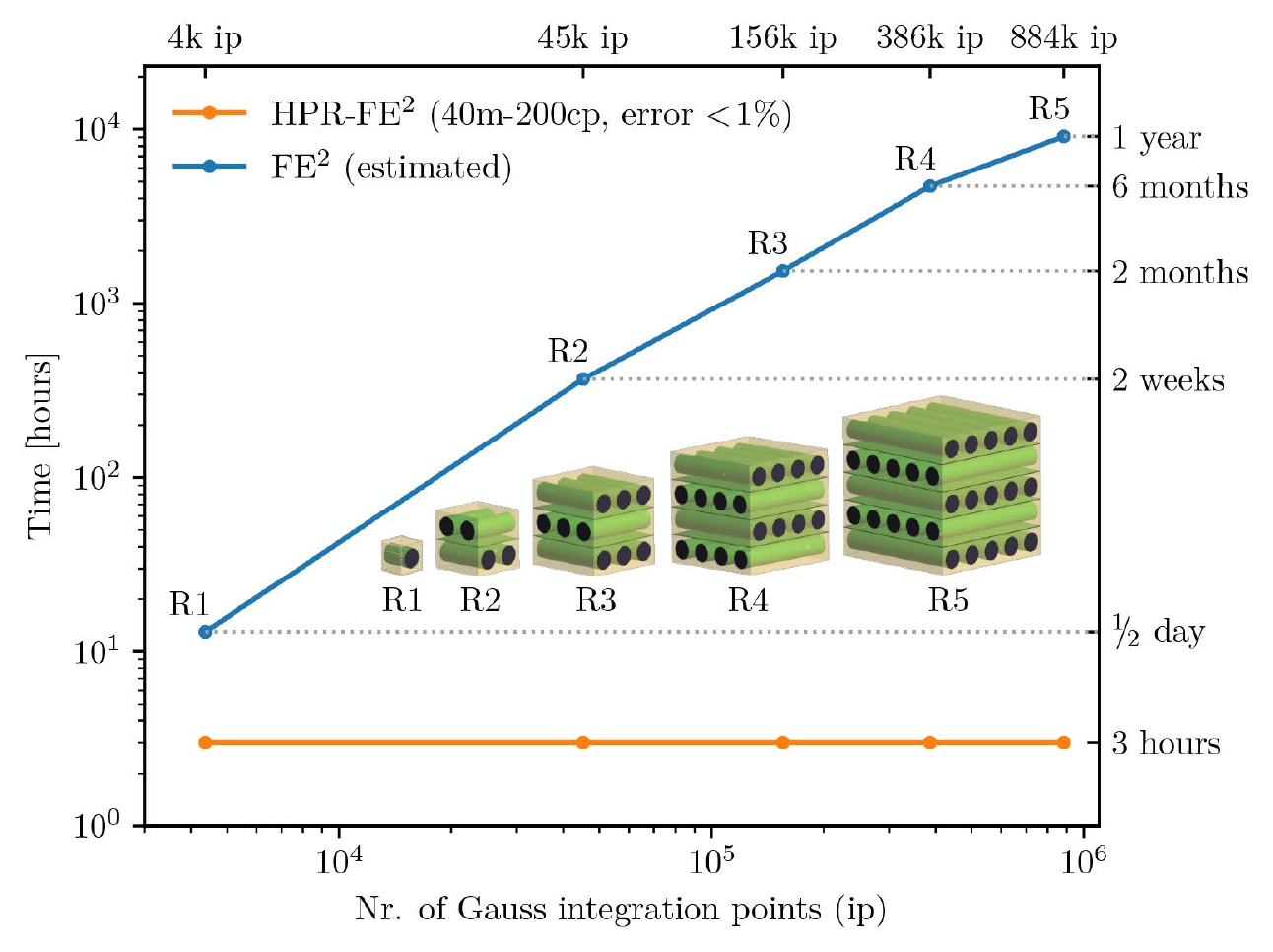}
	\caption{Time for a sequential (one computing thread) \gls{fe2} and \gls{hprfe2} analysis of the coupon and  microcells R1 to R5. The right hand side values are the ``calendar" time, for an intuitive comparison of the time.}
	\label{fig:speedup_time}
\end{figure}

The sequential computational time for a multiscale analysis of the coupon is reported in Figure~\ref{fig:speedup_time}, for the estimated time of \gls{fe2} and the measured time of the \gls{hprfe2} analyses.

For the \gls{fe2} cases, time is estimated as the product of the mean iteration time, the mean number of iterations per step, the number of steps and the total number of Gauss integration points of the coupon discretization.
The \gls{fe2} analysis corresponding to the most complex microcell would demand a year of computation.

The measured time corresponding to the \gls{hprfe2} analyses is close to three hours.
It should be noted that all reduced bases employed for microcells R1 to R5 involve 40 strain modes and 200 cubature points, yielding relative errors below 1\% in the 5 cases.
The error analysis is performed in independent analyses considering only one macroscale point.

Consequently, based on this simple multiscale analysis, it is concluded that the \gls{hprfe2} presented in this work clearly breaks the barrier imposed by the multiplicative cost within hierarchical multiscale analysis, and represents a real chance to export the \gls{fe2} technology to the simulation industry.

%%%%%%%%%%%%%%%%%%%

%%%%%%%%%%%%%%%%%%%
% Section
\section{Conclusions}
The \gls{hprfe2} technology is assessed in this contribution in terms of its benefits for multiscale solutions for a number of cases of industrial relevance, concerning 3D reinforced composite laminate materials.
After an exhaustive evaluation of the presented technique on the above mentioned materials a number of key advantage for the material design industry arise, in comparison to alternative approaches.
\begin{itemize}
	\item \textbf{Mechanical coherency}. The formulation of the multiscale reduced order modelling in terms of the strain field fluctuation turns to be ideal in the context of a \gls{fe2}-like approach, in which downscaling and upscaling contain information of the macro-stains and stresses.
    The strain bases contain only information relevant to the mechanics inheriting the type of boundary conditions employed during the sampling stage, and that is compatible with the strain field, and neglecting the rigid body modes which do not contribute to the constitutive relation.
    Moreover, the strain snapshots are directly Gauss points quantities and not nodal values, which simplifies the formulation and avoids interpolating quantities to the sampling points needed to integrate the strain modes.
    Another important feature of the \gls{hprfe2} relies in the methodology adopted to integrate the strain modes.
    The \gls{roec} integration technique arises from approximating the free energy field (fundamental variational principle) and not its derivatives, as proposed in~\cite{J.A.Hernandez2016}.
    As a result, the problem goal, i.e., find the minimum of the variational principle (free energy), is the actual key of the problem solution unlike in alternative methods based on finding null values of the functional derivatives.
    Consequently, it results in a low-cost and very accurate formulation, leading to the high performance reduction technique \gls{hprfe2}.

    \item \textbf{Numerical consistency}. The accuracy of the approximation, as compared to the \gls{hf} analysis, improves with the number of the reduced strain modes considered in the construction of the strain basis.
    Additionally, an increase of integration points (or energy modes), involved in the integration of the strain modes, provides an improvement of the accuracy with respect to the \gls{hf} solution and, in the limit, tends to the solution provided with a \gls{rom} model including all quadrature points of the original \gls{fe} discretization.
    Noteworthy accuracy is found by adopting a low number of strain modes and cubature points.

    \item \textbf{Model completeness}. Since the material parameters employed within the integration of the reduced strain basis are obtained from the evaluation of the constitutive relation, the link with the physics and the consistent thermodynamics is never lost.
    For this reason, the reduced model can capture well those loading situations (and even thermodynamic irreversible phenomena, like irreversible loading-unloading processes) that were not included during the sampling stage, without the need for including extra modes, as it would happen in other methodologies (e.g., neural network surrogates).

    \item \textbf{Material customization}. The reduced basis can be seen as a discretization support for the integration of the solution.
    In this view, one can keep the integration support and substantially change the model parameters for material design purposes.
    This has been shown in the above sections for significant changes of the material parameters without compromising the quality of the solution.
    Indeed, this can not be taken as a general recipe and radical changes  of the model parameters may, in some cases, demand a more complete strain basis to properly capture the correct solution.
    However, for the type of materials adopted in this study, the material parameters could be remarkably varied without noticing a clear deviation from the \gls{hf} solution.
    The possibility of a reduced model supporting a certain degree of material customizations is considered an overriding quality when the reduced model has to be applied to the design of new materials.

	\item \textbf{Microstructure monitoring and design}. The main feature of \gls{fe2} techniques---modeling the coupled physics in two separated scales---is preserved in the presented approach.
    It relies in monitoring the microstructural behaviour at different key macrostructural points to determine the performance of different phases submitted to stress concentrations and complex loading.
    This is performed by reconstructing the displacement, stress and internal variable fields as indicated in the formulation of the model.
    On this view,
    an ``overall" picture of the material behavior at the low scale is obtained, and
    relevant decisions concerning the microstructural topology and constituents can be taken in order to improve the macrostructural behavior. Again, a desirable feature for the simulation tools used by the material design industry.

    \item \textbf{High computational performance}. A number of tests have been conducted on an engineering coupon with microstructures of increasing  increasing size and discretization features.
    In the presented analyses, the speedup of the approach grows with increasing complexity of the analyzed microstructure, and can reach up to \num{e4} for microstructures with \num{e6} integration points.
    For the case of fully-periodic microcells treated with periodic boundary conditions, this result is not really meaningful since volumes are increased following the periodicity in the three spatial directions and, therefore, essentially contain the same mechanical information.
    However, for the case of minimal kinematic boundary conditions and a sequence of microcells which are not necessarily periodic, a similar observation is done in terms of speedup.
    It is concluded that for the type of materials adopted in the present study, relevant speedup values can be obtained, which allow computations that would typically span several years to be resolved in a few hours and, therefore, being affordable for industrial purposes.
\end{itemize}

To the authors' knowledge, the presented technology represents  a novel contribution towards realistic, accurate and affordable industrial multiscale simulations.
It promises to advance the current standards for material simulation due to a tangible reduction in experimental testing, more complex material behavior that can be modeled with affordable and reliable numerical simulations, seamlessly integrated in current industrial workflows.

%%%%%%%%%%%%%%%%%%%

%%%%%%%%%%%%%%%%%%%
\section*{Acknowledgements}
The authors acknowledge financial support from the Spanish Ministry of Economy and Competitiveness,
through the “Severo Ochoa Programme for Centres of Excellence in R\&D” (CEX2018-000797-S)
and the research grant DPI2017-85521-P for the project “Computational design of Acoustic and Mechanical Metamaterials” (METAMAT).
This research has also received funding from the European Research Council (ERC) under the European Union’s Horizon 2020 research and innovation program (Proof of Concept Grant agreement 874481) through the project “Computational design and prototyping of acoustic metamaterials for target ambient noise reduction” (METACOUSTIC).
The authors also acknowledge the guidance and assistance with the microcells meshes from Dr. Pedro Camanho and Dr. Fermín Otero from INEGI (Portugal) during the preparation of this manuscript.
%%%%%%%%%%%%%%%%%%%

%%%%%%%%%%%%%%%%%%%
\bibliography{ReferenceROM}
%%%%%%%%%%%%%%%%%%%

%%%%%%%%%%%%%%%%%%%
% Section
\appendix
\section{Damage model}
\label{sect:appendix}
A damage model with an elastic domain accounting for traction-only and a bilinear hardening law is adopted for describing the constitutive response of the components at the microscale.
The main features of the model are described in Table \ref{Tab_Damage}.
The employed constitutive model is extensively described in~\cite{oliver2015continuum} for the case of linear softening behaviour.
However, a short explanation of the parameters and symbols used in the table is included in the following lines: $\psi$ is the free energy and $\boldsymbol{C}$ is the Hooke's elasticity tensor written in terms of the Lamé parameters $\lambda$ and $\mu$, with $\mathds{I}$ and $\mathds{1}$  the forth- and second-order identity tensor.
The isotropic damage variable $d$  depends on the internal variable $r$.
This internal variable determines the size of the elastic domain through the damage function $g(\boldsymbol{\varepsilon},r)$.
The term $\tau_{\varepsilon}$ in the definition of $g$ is a generalized norm of the strains $\boldsymbol{\varepsilon}$ accounting for the positive effective principal stresses $\bar{\boldsymbol{\sigma}}^+$.
A bilinear hardening law $q(r)$ is considered in the present study with two values  of the hardening parameter $H$.
It is worth mentioning that no regularization is required since softening behaviour is not included in the present study.
\renewcommand{\arraystretch}{1.5}
\begin{table}[htp]
    \begin{center}
        \begin{tabular} {l  l }
            \toprule
            Free energy  $\psi$: &
            $ \psi(\boldsymbol{\varepsilon} ,d) = (1-d) \psi_0;
            \quad
            \psi_0 = \frac{1}{2} \boldsymbol{\varepsilon} : \boldsymbol{C} : \boldsymbol{\varepsilon}
            \label{free_energy}$
            \\
            &
            $\boldsymbol{C} = 2\mu \mathds{I} + \lambda(\mathds{1} \otimes \mathds{1} )
            \label{hooke_tensor} $
            \\
            \midrule
            Damage variable $d$:  &
            $ d(r) =1 - \frac{q(r)}{r}; \quad q \geq 0; \quad r\geq 0
            \label{damage_variable} $
            \\
            \midrule
            Constitutive equation:  &
            $\boldsymbol{\sigma}
            =
            (1-d)\boldsymbol{C} : \boldsymbol{\varepsilon}
            =
            \frac{q}{r}\underbrace{\boldsymbol{C} : \boldsymbol{\varepsilon}}_{\bar{\boldsymbol{\sigma}}}
            =
            \frac{q}{r}{\bar{\boldsymbol{\sigma}}}
            \label{constitutive_equation_l1} $
            \\
            \midrule
            Damage function $g$: &
            $g(\boldsymbol{\varepsilon}, r)
            =
            \tau_{\varepsilon}(\boldsymbol{\varepsilon}) - r;
            \quad
            \tau_{\varepsilon}(\boldsymbol{\varepsilon})
            =
            \sqrt{ \bar{\boldsymbol{\sigma}}^+ : \boldsymbol{\varepsilon}}
            \label{damage_function_l1} $
            \\
            \midrule
            Initial conditions: &
            $r_0
            =
            \frac{\sigma_{\text{e}}}{\sqrt{E}};
            \quad
            q_0 = r_0
            \label{internal_varable} $
            \\
            \midrule
            Loading-unloading conditions: &
            $\dot{r} \geq 0;
            \quad
            g \leq 0;
            \quad
            \dot{r} g = 0
            \label{loading_condition_l1} $
            \\
            \midrule
            Hardening law $q$:&
            $\dot{q} = H(r) \dot{r};
            \quad
            H > 0
            \quad
            \text{(as defined in Figure~\ref{fig_hardening})}
            \label{softening_law_l1} $
            \\
            \bottomrule
        \end{tabular}
        \caption{Material damage model adopted for the microscale description.}
        \label{Tab_Damage}
    \end{center}
\end{table}
\begin{figure}[htbp]
	\centering
	\includegraphics[width=0.7\linewidth]{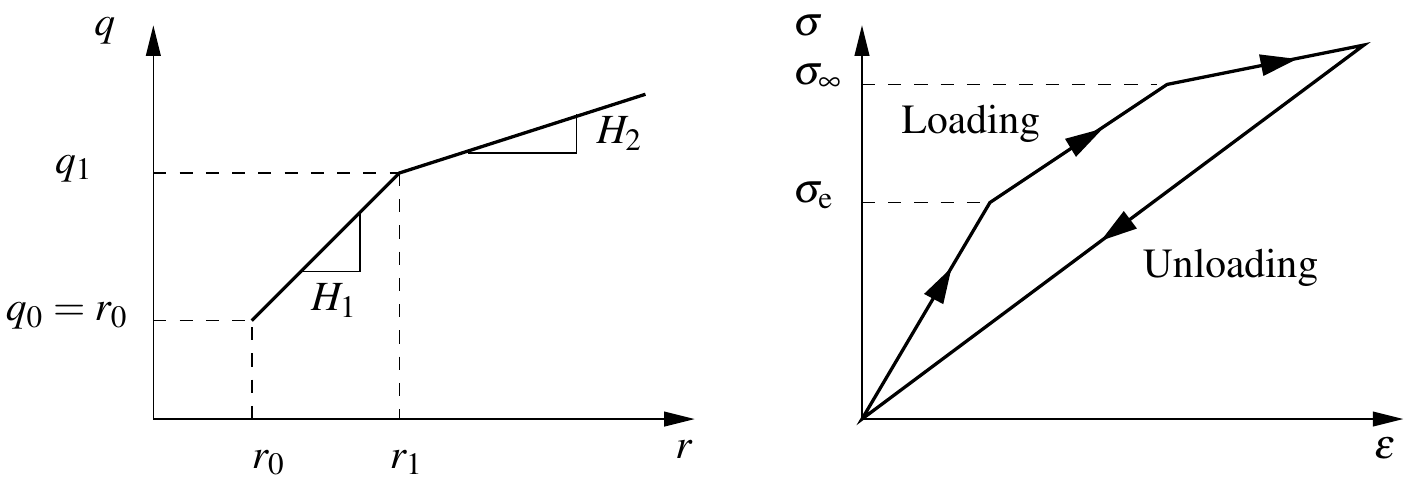}
	\caption{Bilinear hardening law (left). Loading-unloading response in uniaxial tension test (right).}
	\label{fig_hardening}
\end{figure}
%
%%%%%%%%%%%%%%%%%%%

\end{document}